\documentclass[journal]{IEEEtran}
\usepackage{cite}
\ifCLASSINFOpdf
  \usepackage[pdftex]{graphicx}
  \usepackage{epstopdf}
  \usepackage{subfigure}
  \graphicspath{{./}{./figures/}}
\else
\fi
\usepackage{xcolor}

\usepackage{amsmath}
\usepackage{amssymb}
\allowdisplaybreaks
\newcommand*{\qed}{\hfill\ensuremath{\square}}%
\newcommand*\diff{\mathop{}\!\mathrm{d}}%

\newtheorem{theorem}{Theorem}[section]

\newtheorem{lemma}[theorem]{Lemma}

\newtheorem{assum}[theorem]{Assumption}
\newtheorem{remark}[theorem]{Remark}

\begin{document}
\title{Feedback Stabilization of a Class of Diagonal Infinite-Dimensional Systems with Delay Boundary Control}

\author{Hugo~Lhachemi and Christophe~Prieur% <-this % stops a space
\thanks{Hugo Lhachemi is with the School of Electrical and Electronic Engineering, University College Dublin, Dublin, Ireland. e-mail: hugo.lhachemi@ucd.ie. 
Christophe Prieur is with Univ. Grenoble Alpes, CNRS, Grenoble-INP, GIPSA-lab, F-38000, Grenoble, France. email: christophe.prieur@gipsa-lab.fr.
This publication has emanated from research supported in part by a research grant from Science Foundation Ireland (SFI) under grant number 16/RC/3872 and is co-funded under the European Regional Development Fund and by I-Form industry partners.}% <-this % stops a space
}

%\author{Hugo~Lhachemi% <-this % stops a space
%\thanks{Hugo Lhachemi is with the School of Electrical and Electronic Engineering, University College Dublin, Dublin, Ireland. e-mail: hugo.lhachemi@ucd.ie.}% <-this % stops a space
%}


\markboth{Preprint}%
{Lhachemi \MakeLowercase{\textit{et al.}}: Feedback Stabilization of a Class of Diagonal Infinite-Dimensional Systems with Delay Boundary Control}



\maketitle

\begin{abstract}
This paper studies the boundary feedback stabilization of a class of diagonal infinite-dimensional boundary control systems. In the studied setting, the boundary control input is subject to a constant delay while the open loop system might exhibit a finite number of unstable modes. The proposed control design strategy consists in two main steps. First, a finite-dimensional subsystem is obtained by truncation of the original Infinite-Dimensional System (IDS) via modal decomposition. It includes the unstable components of the infinite-dimensional system and allows the design of a finite-dimensional delay controller by means of the Artstein transformation and the pole-shifting theorem. Second, it is shown via the selection of an adequate Lyapunov function that 1) the finite-dimensional delay controller successfully stabilizes the original infinite-dimensional system ; 2) the closed-loop system is exponentially Input-to-State Stable (ISS) with respect to distributed disturbances. Finally, the obtained ISS property is used to derive a small gain condition ensuring the stability of an IDS-ODE interconnection.
\end{abstract}


% \fi
\IEEEpeerreviewmaketitle

\section{Introduction}\label{sec: Introduction}

Feedback control of finite-dimensional systems in the presence of input delays has been extensively investigated~\cite{artstein1982linear,richard2003time}. The extension of this topic to Infinite-Dimensional Systems (IDS), and in particular to Partial Differential Equations (PDEs), has attracted much attention in the recent years.

There exist essentially two types of control inputs for infinite-dimensional systems: bounded and unbounded control operators. The stability of linear and semilinear infinite-dimensional system under time-varying delayed feedback acting via a bounded linear control operator has been studied, e.g., in~\cite{fridman2009exponential,solomon2015stability}. In this paper, we are interested in the second type of control, i.e., when the control input acts on the system via an unbounded operator. For PDEs, such a setting takes the form of a control acting in the boundary conditions.

Unbounded control operators have been considered in the stability study of various PDEs. The cases of the heat~\cite{nicaise2009stability} and wave~\cite{nicaise2008stabilization,nicaise2007stabilization,nicaise2009stability} equations were studied via Lyapunov methods for slow time vaying delays. The cases of a parabolic PDE and a second-order evolution equation were reported in~\cite{wang2018delay} and~\cite{fridman2010stabilization}, respectively. The extension to a delayed ODE–heat cascade under actuator saturation was reported in~\cite{kang2017boundary}.

In this paper, we are interested in the boundary feedback stabilization of a class of diagonal infinite dimensional boundary control systems in the presence of a constant input delay. Specifically, we consider the case of a boundary control system~\cite{Curtain2012} for which the associated disturbance free operator is a Riesz-spectral operator admitting a finite number of unstable eigenvalues. The control design objective consists in the feedback stabilization of the system by means of a delay boundary control. 

One of the very first contributions on input delayed unstable PDEs deals with a reaction-diffusion equation~\cite{krstic2009control} where the controller was designed by resorting to the backstepping technique. The approach adopted in this paper differs. It relies on the following three steps procedure initially reported in~\cite{russell1978controllability}: 1) obtaining a finite-dimensional subsystem capturing the unstable modes by truncation of the original infinite-dimensional via a modal decomposition ; 2) design of a finite dimensional control law that stabilizes the finite-dimensional unstable part of the system ; 3) use of an adequate Lyapunov function to assess that the designed control law stabilizes the original infinite-dimensional system. Such a control design strategy was successfully applied to the stabilization of semilinear heat~\cite{coron2004global} and wave~\cite{coron2006global} equations via (undelayed) boundary feedback control. The extension of this design procedure to the delay feedback control of a linear reaction-diffusion equation was reported in~\cite{prieur2018feedback,tanwani2018local}. The delayed finite dimensional model was obtained via spectral reduction. Then, the control law was computed by applying the Artstein transformation~\cite{artstein1982linear,richard2003time} and by resorting to the classical pole-shifting theorem. A distinguished feature is that, under the knowledge of the constant delay $D \geq 0$, the obtained finite-dimensional control law amounts stabilizing the closed-loop system, whatever the value of the time-delay $D$ may be.

In this context, the contribution of the present paper is fourfold.
\begin{enumerate}
\item We generalize the approach developed in~\cite{prieur2018feedback} for the delay feedback control of a linear reaction-diffusion equation with one-dimensional control input to the general case of the delay boundary feedback stabilization of a class of diagonal infinite dimensional boundary control systems with finitely many unstable modes and finite dimensional input. The control design strategy relies on the design of the feedback control law based on a finite-dimensional truncated part of the original system. The truncation is performed via a spectral decomposition used to capture the unstable modes of the system. The control law is then obtained based on this finite-dimensional subsystem with delay control input by means of the Artstein transformation and the pole-shifting theorem. The exponential stability of the resulting closed-loop infinite dimensional system is assessed via the introduction of a suitable Lyapunov function.
\item In~\cite{coron2004global,coron2006global,prieur2018feedback} the control design was performed on the time derivative $v = \dot{u}$ of the actual input signal $u$. Thus, the application of the control law required an \emph{a posteriori} integration of $v$ to obtain the actual control input $u$. In this paper, we propose a simplification of the control law that avoids such an \emph{a posteriori} integration. Such a simplification is allowed by an adequate spectral decomposition that only involves the value of the control input while avoiding the occurrence of its time derivative.
\item We show that the resulting closed-loop system is exponentially Input-to-State Stable (ISS)\cite{sontag1989smooth} with respect to distributed disturbances acting via a bounded operator. 
\item Taking advantage of the ISS property of the closed-loop infinite-dimensional system, we derive a small gain condition ensuring the stability of an IDS-ODE interconnection. We follow here the methodology presented in~\cite{karafyllis2018input} that relies on the conversion of the ISS estimates satisfied by each component of the interconnection into fading memory estimates\cite[Lemma~7.1]{karafyllis2018input}. However, such a conversion does not apply to the studied closed-loop infinite-dimensional system due to the time-varying nature of the control strategy. This pitfall is avoided by working directly with the Lyapunov function instead of the trajectories of the system.
\end{enumerate}

The remainder of this paper is organized as follows. Both problem setting and control objectives are introduced in Section~\ref{sec: Problem setting and control objective}. The comprehensive construction of the control strategy is presented in Section~\ref{sec: construction of the feedback control strategy}. It consists first in the spectral decomposition of the problem in order to obtain a finite-dimensional model capturing the unstable modes (Subsection~\ref{subsec: spectral decomposition}) and then the design of the finite dimensional controller stabilizing the obtained truncated subsystem. The study of the ISS property of the resulting closed-loop infinite-dimensional system is carried out via the introduction of an adequate Lyapunov function in Section~\ref{sec: Study of the infinite-dimensional closed-loop system}. We take advantage of these results to derive in Section~\ref{sec: Application to the stability analysis of a closed-loop interconected IDS-ODE system} a small gain condition ensuring the stability of an IDS-ODE interconnection. In Section~\ref{sec: case study}, we check the assumptions on a IDS-ODE system and in particular the small gain condition. The obtained numerical results are compliant with the theoretical predictions. Finally, concluding remarks are provided in Section~\ref{sec: conclusion}.

\section{Problem setting and control objective}\label{sec: Problem setting and control objective}
Throughout the paper, we assume that $(\mathcal{H},\left< \cdot , \cdot \right>_\mathcal{H})$ is a separable Hilbert space over the field $\mathbb{K}$, which is either $\mathbb{R}$ or $\mathbb{C}$. All the finite-dimensional spaces $\mathbb{K}^p$ are endowed with the usual euclidean inner product $\left<x,y\right> = x^* y$ and the associated 2-norm $\Vert x \Vert = \sqrt{\left<x,x\right>} = \sqrt{x^* x}$, where $x^* = \overline{x}^\top$. For any matrix $M \in \mathbb{K}^{p \times q}$, $\Vert M \Vert$ stands for the induced norm of $M$ associated with the above 2-norms.

\subsection{Problem setting}
We consider the abstract boundary control systems~\cite{Curtain2012} with delayed boundary control 
\begin{equation}\label{def: boundary control system}
\left\{\begin{split}
\dfrac{\mathrm{d} X}{\mathrm{d} t}(t) & = \mathcal{A} X(t) + d(t) , & t \geq 0 \\
\mathcal{B} X (t) & = u_D(t) \triangleq u(t-D) , & t \geq 0 \\
X(0) & = X_0 
\end{split}\right.
\end{equation}
with
\begin{itemize}
\item $\mathcal{A} : D(\mathcal{A}) \subset \mathcal{H} \rightarrow \mathcal{H}$ a linear (unbounded) operator;
\item $\mathcal{B} : D(\mathcal{B}) \subset \mathcal{H} \rightarrow \mathbb{K}^m$ with $D(\mathcal{A}) \subset D(\mathcal{B})$ a linear boundary operator;
\item $d : \mathbb{R}_+ \rightarrow \mathcal{H}$ a distributed disturbance;
\item $u : [-D,+\infty) \rightarrow \mathbb{K}^m$, with a known constant delay $D >0$ and $\left. u \right\vert_{[-D,0)} = 0$, the boundary control.
\end{itemize}
We assume that $(\mathcal{A},\mathcal{B})$ is a boundary control system, i.e.,
\begin{enumerate}
\item the disturbance free operator $\mathcal{A}_0$, defined over the domain $D(\mathcal{A}_0) \triangleq D(\mathcal{A}) \cap \mathrm{ker}(\mathcal{B})$ by $\mathcal{A}_0 \triangleq \left.\mathcal{A}\right|_{D(\mathcal{A}_0)}$, is the generator of a $C_0$-semigroup $S$ on $\mathcal{H}$;
\item there exists a bounded operator $B \in \mathcal{L}(\mathbb{K}^m,\mathcal{H})$, called a lifting operator, such that $\mathrm{R}(B) \subset D(\mathcal{A})$, $\mathcal{A}B \in \mathcal{L}(\mathbb{K}^m,\mathcal{H})$, and $\mathcal{B}B = I_{\mathbb{K}^m}$.
\end{enumerate}
It is recalled that $\mathrm{ker}(\mathcal{B})$ is the kernel of $\mathcal{B}$ while $\mathrm{R}(B)$ stands for the range of $B$. We make the following assumptions.

\begin{assum}\label{assum: A1}
The disturbance free operator $\mathcal{A}_0$ is a Riesz spectral operator~\cite{Curtain2012}, i.e., is a linear and closed operator with simple eigenvalues $\lambda_n$ and corresponding eigenvectors $\phi_n \in D(\mathcal{A}_0)$, $n \in \mathbb{N}^*$, that satisfy:
\begin{enumerate}
\item $\left\{ \phi_n , \; n \in \mathbb{N}^* \right\}$ is a Riesz basis~\cite{christensen2016introduction}:
\begin{enumerate}
\item $\overline{ \underset{n\in\mathbb{N}^*}{\mathrm{span}_\mathbb{K}} \;\phi_n } = \mathcal{H}$;
\item there exist constants $m_R, M_R \in \mathbb{R}_+^*$ such that for all $N \in \mathbb{N}^*$ and all $\alpha_1 , \ldots , \alpha_N \in \mathbb{K}$,
\begin{equation}\label{eq: Riesz basis - inequality}
m_R \sum\limits_{n=1}^{N} \vert \alpha_n \vert^2
\leq
\left\Vert \sum\limits_{n=1}^{N} \alpha_n \phi_n \right\Vert_\mathcal{H}^2
\leq
M_R \sum\limits_{n=1}^{N} \vert \alpha_n \vert^2 .
\end{equation}
\end{enumerate}
\item The closure of $\{ \lambda_n , \; n \in \mathbb{N}^* \}$ is totally disconnected, i.e. for any distinct $a,b \in \overline{ \{ \lambda_n , \; n \in \mathbb{N}^* \} }$, $[a,b] \not\subset \overline{ \{ \lambda_n , \; n \in \mathbb{N}^* \} }$.
\end{enumerate}
\end{assum}

\begin{assum}\label{assum: A2}
There exist $N_0 \in \mathbb{N}^*$ and $\alpha \in \mathbb{R}_+^*$ such that $\operatorname{Re} \lambda_n \leq - \alpha$ for all $n \geq N_0 + 1$. 
\end{assum}

\begin{remark}
Note that Assumption~\ref{assum: A2} is equivalent to:
\begin{itemize}
\item the number of unstable eigenvalues is finite, i.e., $\mathrm{Card} (\{ \lambda_n \,:\, \operatorname{Re} \lambda_n \geq 0\}) < \infty$ ;
\item the set composed of the real part of the stable eigenvalues is not accumulating at 0, i.e., $\sup\limits_{\operatorname{Re}\lambda_n < 0} \operatorname{Re} \lambda_n < 0$. 
\end{itemize}
\end{remark}

From the well-known properties of the Riesz-basis (see, e.g.,~\cite{christensen2016introduction}), we introduce $\left\{ \psi_n , \; n \in \mathbb{N}^* \right\}$ the biorthogonal sequence associated with the Riesz basis $\left\{ \phi_n , \; n \in \mathbb{N}^* \right\}$, i.e., $\left< \phi_k , \psi_l \right>_\mathcal{H} = \delta_{k,l}$. Then, the following series expansion holds true.
\begin{equation*}
\forall x \in \mathcal{H}, \qquad
x = \sum\limits_{n \geq 1} \left< x , \psi_n \right>_\mathcal{H} \phi_n .
\end{equation*}
Furthermore, as $\mathcal{A}_0$ is assumed to be a Riesz-spectral operator, then $\psi_n$ is an eigenvector of the adjoint operator $\mathcal{A}_0^*$ associated with the eigenvalue $\overline{\lambda_n}$. 

\subsection{Control objective}
The control objective is twofold. First, in the absence of distributed disturbance (i.e., $d=0$), the control objective is to design a control law $u$ that exponentially stabilizes (if at least one eigenvalue has a non negative real part) and modify the pole placement associated with $\lambda_1,\ldots,\lambda_{N_0}$ for (\ref{def: boundary control system}). Second, the control law must ensure the ISS property of the closed-loop system with respect to the distributed disturbance $d$. 

Because we are only concerned in controlling the system from the starting time $t = 0$, we assume that the system is uncontrolled for $t < 0$. This is why it is imposed $\left. u \right\vert_{[-D,0)} = 0$. Therefore, due to the delay $D$ in the control input of (\ref{def: boundary control system}), the system remains open-loop for $t < D$ while the effect of the control input will have an impact on the system only at times $t \geq D$.

Note that the $N_0 \in \mathbb{N}^*$ and $\alpha > 0$ provided by Assumption~\ref{assum: A2} are not unique. For instance, one could select $N_0 \in \mathbb{N}^*$ such that $\lambda_1,\ldots,\lambda_{N_0}$ are all with non negative real part. In this case, the control design reduces to stabilize the unstable part of the system. Nevertheless, one could also want to improve the decay rate or the damping of certain of the stable open-loop eigenvalues. In this case, $\lambda_1,\ldots,\lambda_{N_0}$ would include all the unstable eigenvalues and certain selected stable eigenvalues of the open-loop system.

\section{Construction of the feedback control strategy}\label{sec: construction of the feedback control strategy}
In order to derive the control law, we make in this section the \emph{a priori} assumption that $u \in \mathcal{C}^2([-D,+\infty);\mathbb{K}^m)$. This assumption is necessary to ensure the existence of classical solutions of (\ref{def: boundary control system}), and thus to proceed to the upcoming computations (see, e.g.,~\cite{Curtain2012}). Therefore, the construction of the control law must ensure that such a regularity property holds true. This will be assessed in the next section.

\subsection{Spectral decomposition}\label{subsec: spectral decomposition}
Assuming that  $u_D \in \mathcal{C}^2([0,+\infty);\mathbb{K}^m)$, $X_0 \in D(\mathcal{A})$ such that $\mathcal{B}X_0 = u_D(0) = 0$ (i.e., $X_0 \in D(\mathcal{A}_0)$), and $d \in \mathcal{C}^1(\mathbb{R}_+;\mathcal{H})$, we denote by $X \in \mathcal{C}^0(\mathbb{R}_+;D(\mathcal{A})) \cap \mathcal{C}^1(\mathbb{R}_+;\mathcal{H})$ the unique classical solution of (\ref{def: boundary control system}). Then, we introduce $c_n(t) \triangleq \left< X(t) , \psi_n \right>_\mathcal{H}$ the projection of $X(t)$ into the Riezs basis $\left\{ \phi_n , \; n \in \mathbb{N}^* \right\}$, i.e. (see~\cite{christensen2016introduction}),
\begin{equation}\label{eq: proj Riesz basis}
X(t) 
= \sum\limits_{n \in \mathbb{N}^*} \left< X(t) , \psi_n \right>_\mathcal{H} \phi_n 
= \sum\limits_{n \in \mathbb{N}^*} c_n(t) \phi_n 
.
\end{equation}
We also introduce $d_n(t) \triangleq \left< d(t) , \psi_n \right>_\mathcal{H}$. Then $c_n \in \mathcal{C}^1(\mathbb{R}_+;\mathbb{K})$ and, following~\cite{lhachemi2018iss}, we infer from (\ref{def: boundary control system}) that, for all $t \geq 0$,
\begin{align}
& \dot{c}_n(t) \nonumber \\
& = \left< \dfrac{\mathrm{d}X}{\mathrm{d}t}(t) , \psi_n \right>_\mathcal{H} \nonumber \\
& = \left< \mathcal{A} X(t) , \psi_n \right>_\mathcal{H} + \left< d(t) , \psi_n \right>_\mathcal{H} \nonumber \\
& = \left< \mathcal{A} \left\{ X(t) - B u_D(t) \right\} , \psi_n \right>_\mathcal{H} + \left< \mathcal{A} B u_D(t) , \psi_n \right>_\mathcal{H} + d_n(t) \nonumber \\
& = \left< \mathcal{A}_0 \left\{ X(t) - B u_D(t) \right\} , \psi_n \right>_\mathcal{H} + \left< \mathcal{A} B u_D(t) , \psi_n \right>_\mathcal{H} + d_n(t) \nonumber \\
& = \left< X(t) - B u_D(t) , \mathcal{A}_0^* \psi_n \right>_\mathcal{H} + \left< \mathcal{A} B u_D(t) , \psi_n \right>_\mathcal{H} + d_n(t) \nonumber \\
& = \left< X(t) - B u_D(t) , \overline{\lambda_n} \psi_n \right>_\mathcal{H} + \left< \mathcal{A} B u_D(t) , \psi_n \right>_\mathcal{H} + d_n(t) \nonumber \\
& = \lambda_n c_n(t) - \lambda_n \left< B u_D(t) , \psi_n \right>_\mathcal{H} + \left< \mathcal{A} B u_D(t) , \psi_n \right>_\mathcal{H} + d_n(t) , \label{eq: coeff in Riesz basis ODE}
\end{align}
where it has been used that $\mathcal{B} \left\{ X(t) - B u_D(t) \right\} = u_D(t) - u_D(t) = 0$, showing that $X(t) - B u_D(t) \in D(\mathcal{A}) \cap \mathrm{ker}(\mathcal{B}) = D(\mathcal{A}_0)$.

\begin{remark}
It is interesting to note that the ODE (\ref{eq: coeff in Riesz basis ODE}) describing the time evolution of the coefficient $c_n(t) = \left< X(t) , \psi_n \right>_\mathcal{H}$ only involves the delayed control input $u_D(t)$ while avoiding the occurrence of its time derivative $\dot{u}_D(t)$. Therefore, whereas it was necessary in~\cite{coron2004global,coron2006global,prieur2018feedback}, due to the presence of the term $\dot{u}_D(t)$ in the ODEs resulting from the spectral decomposition, to augment the state of the finite-dimensional subsystem and to use $\dot{u}_D(t)$ as a control input, we avoid here such a procedure. This yields a simplification of the control law by avoiding an \emph{a posteriori} integration of $\dot{u}$ to obtain the actual control law $u$.
\end{remark}

Let $\mathcal{E} = (e_1,e_2,\ldots,e_m)$ be the canonical basis of $\mathbb{K}^m$, and consider the projections $u_1,u_2,\ldots,u_m \in \mathcal{C}^2([-D,+\infty);\mathbb{K})$ such that
\begin{equation*}
u 
= \sum\limits_{k=1}^{m} u_k e_k 
= \begin{bmatrix}
u_1 \\ \vdots \\ u_{m}
\end{bmatrix} .
\end{equation*}
Introducing $b_{n,k} \triangleq - \lambda_n \left< B e_k , \psi_n \right>_\mathcal{H} + \left< \mathcal{A} B e_k , \psi_n \right>_\mathcal{H}$, we obtain from (\ref{eq: coeff in Riesz basis ODE}) that
\begin{align*}
\dot{c}_n(t) 
& = \lambda_n c_n(t) + \sum\limits_{k=1}^m b_{n,k} u_{D,k}(t)  + \left< d(t) , \psi_n \right>_\mathcal{H} .
\end{align*}
Then, the following linear ODE with delay input holds true for all $t \geq 0$
\begin{equation}\label{eq: ODE satisfies by Y}
\dot{Y}(t) 
= A_{N_0} Y(t) + B_{N_0} u_D(t) + D_{N_0}(t) ,
\end{equation}
where $A_{N_0} = \mathrm{diag}(\lambda_1,\ldots,\lambda_{N_0}) \in \mathbb{K}^{N_0 \times N_0}$, $B_{N_0} = (b_{n,k})_{1 \leq n \leq N_0 , 1 \leq k \leq m} \in \mathbb{K}^{N_0 \times m}$,
\begin{equation*}
Y(t) = 
\begin{bmatrix}
c_1(t) \\ \vdots \\ c_{N_0}(t)
\end{bmatrix}
=
\begin{bmatrix}
\left< X(t) , \psi_1 \right>_\mathcal{H} \\ \vdots \\ \left< X(t) , \psi_{N_0} \right>_\mathcal{H}
\end{bmatrix} 
\in \mathbb{K}^{N_0} ,
\end{equation*}
and
\begin{equation}\label{eq: def DN0}
D_{N_0}(t) = 
\begin{bmatrix}
d_1(t) \\ \vdots \\ d_{N_0}(t)
\end{bmatrix}
=
\begin{bmatrix}
\left< d(t) , \psi_1 \right>_\mathcal{H} \\ \vdots \\ \left< d(t) , \psi_{N_0} \right>_\mathcal{H}
\end{bmatrix} 
\in \mathbb{K}^{N_0} .
\end{equation}
Note that the norm of $D_{N_0}(t)$ can be bounded above in function of the norm of the full distributed disturbance $d(t)$ as follows. For all $t \geq 0$, we have
\begin{align}
\Vert D_{N_0}(t) \Vert^2 
= \sum\limits_{k=1}^{N_0} \vert \left< d(t) , \psi_k \right>_\mathcal{H} \vert^2
& \leq \sum\limits_{k \geq 1} \vert \left< d(t) , \psi_k \right>_\mathcal{H} \vert^2 \nonumber \\
& \overset{(\ref{eq: Riesz basis - inequality})}{\leq} \dfrac{1}{m_R} \Vert d(t) \Vert_\mathcal{H}^2 . \label{eq: above estimate of D_N0}
\end{align}

The finite-dimensional linear ODE (\ref{eq: ODE satisfies by Y}) captures the part of the dynamics of (\ref{def: boundary control system}) that must me stabilized/controlled by the feedback control $u$. The idea consists in first designing a control law that exponentially stabilizes the linear ODE (\ref{eq: ODE satisfies by Y}). Then, we assess that the proposed control law amounts stabilizing the original infinite-dimensional system (\ref{def: boundary control system}) by means of an adequate Lyapunov function.

\subsection{Stabilization of the finite-dimensional subsystem}
At this point, we need to design a control law that stabilizes the linear ODE with input delay (\ref{eq: ODE satisfies by Y}). First, we resort to the Artstein model reduction~\cite{artstein1982linear,richard2003time} to obtain an equivalent linear ODE that is free of delay. Specifically, we introduce for all $t \geq 0$,
\begin{align*}
Z(t) & = Y(t) + \int_{t-D}^{t} e^{(t-s-D)A_{N_0}} B_{N_0} u(s) \diff s \\
     & = Y(t) + \int_{0}^{D} e^{-\tau A_{N_0}} B_{N_0} u(t-D+\tau) \diff \tau .
\end{align*}
Straightforward computations show that we have for all $t \geq 0$,
\begin{equation*}
\dot{Z}(t) = A_{N_0} Z(t) + e^{-D A_{N_0}} B_{N_0} u(t) + D_{N_0}(t) .
\end{equation*}
As $e^{-D A_{N_0}}$ is invertible and commutes with $A_{N_0}$, the pair $(A_{N_0},e^{-D A_{N_0}} B_{N_0})$ satisfies the Kalman condition if and only if the pair $(A_{N_0},B_{N_0})$ satisfies the Kalman condition. Consequently, in order to be able to apply the pole-shifting theorem, we make the following assumption.

\begin{assum}\label{assum: A3}
$(A_{N_0},B_{N_0})$ satisfies the Kalman condition.
\end{assum}

\begin{remark}\label{rem: A3}
In the case of a one-dimensional control input, i.e., $m = 1$, we have that
\begin{align*}
& \mathrm{det}(B_{N_0}, A_{N_0}B_{N_0}, \ldots, A_{N_0}^{N_0-1} B_{N_0} ) \\
& \qquad = \prod\limits_{n=1}^{N_0} b_{n,1} \times \mathrm{VdM}(\lambda_1,\ldots,\lambda_{N_0}) ,
\end{align*} 
where $\mathrm{VdM}(\lambda_1,\ldots,\lambda_{N_0})$ is the Van der Monde determinant associated with $\lambda_1,\ldots,\lambda_{N_0}$. Therefore Assumption~\ref{assum: A3} is fulfilled if and only if $\lambda_1,\ldots,\lambda_{N_0}$ are all distinct and $b_{n,1} \neq 0$ for all $1 \leq n \leq N_0$. In the general case $m \geq 1$, we can easily apply the PBH test~\cite{zhou1998essentials} due to the diagonal nature of the matrix $A_{N_0}$. Assume without loss of generality that $\lambda_1,\ldots,\lambda_{N_0}$ are ordered such that there exist $n_1,\ldots,n_p \in \mathbb{N}^*$ with $n_1 + \ldots + n_p = N_0$ such that 1) for all $1 \leq l \leq p$, $\lambda_{s_{l-1} + 1} = \lambda_{s_{l-1} + 2} = \ldots = \lambda_{s_l}$ ; 2) $l_1 \neq l_2$ implies $\lambda_{s_{l_1}} \neq \lambda_{s_{l_2}}$, where $s_{l} = n_1 + n_2 + \ldots + n_l$. Then, Assumption~\ref{assum: A3} is fulfilled if and only if $\mathrm{rank} [ (b_{n,k})_{s_{l-1} + 1 \leq n \leq s_l , 1 \leq k \leq m} ] = n_l$. In particular, it requires the necessary condition that $n_l \leq m$ for all $1 \leq l \leq p$.
\end{remark}

\begin{remark}
Note that $b_{n,k}$ is computed based on the selection of a given lifting operator $B$. Even if such a lifting operator is not unique, the quantity $b_{n,k}$ is actually independent of the particularly selected lifting operator. Indeed, let $B$ and $\tilde{B}$ be two distinct lifting operators associated with $(\mathcal{A},\mathcal{B})$. Then, introducing $\hat{B} = B - \tilde{B}$, one has $\mathcal{B}\hat{B} = \mathcal{B} B - \mathcal{B} \tilde{B} = I_{\mathbb{K}^m} - I_{\mathbb{K}^m} = 0$. Thus, $R(\hat{B}) \subset D(\mathcal{A}) \cap \mathrm{ker}(\mathcal{B}) = D(\mathcal{A}_0)$ and we obtain that
\begin{align*}
\left< \mathcal{A} \hat{B} e_k , \psi_n \right>_\mathcal{H}
& = \left< \mathcal{A}_0 \hat{B} e_k , \psi_n \right>_\mathcal{H}
= \left< \hat{B} e_k , \mathcal{A}_0^* \psi_n \right>_\mathcal{H} \\
& = \left< \hat{B} e_k , \overline{\lambda_n} \psi_n \right>_\mathcal{H}
= \lambda_n \left< \hat{B} e_k , \psi_n \right>_\mathcal{H} .
\end{align*}
We deduce the claimed result, i.e., 
\begin{align*}
& - \lambda_n \left< B e_k , \psi_n \right>_\mathcal{H} + \left< \mathcal{A} B e_k , \psi_n \right>_\mathcal{H} \\
& \hspace{1cm} = - \lambda_n \left< \tilde{B} e_k , \psi_n \right>_\mathcal{H} + \left< \mathcal{A} \tilde{B} e_k , \psi_n \right>_\mathcal{H} .
\end{align*}
Therefore, the commandability property of the pair $(A_{N_0},B_{N_0})$ is an intrinsic property of the boundary control system $(\mathcal{A},\mathcal{B})$ in the sense that it does not depend on the selection of a particular lifting operator $B$.
\end{remark}

Assuming that Assumption~\ref{assum: A3} holds true, we can find a feedback gain $K \in \mathbb{K}^{m \times N_0}$ and $P \in \mathcal{H}_{N_0}^{+*}$ Hermitian definite positive such that $A_\mathrm{cl} \triangleq A_{N_0} + e^{-D A_{N_0}} B_{N_0} K$ is Hurwitz with desired pole placement and 
\begin{equation*}
A_\mathrm{cl}^* P + P A_\mathrm{cl} = -I_{N_0} .
\end{equation*}
Then, a natural choice for the control input would be $u(t) = \chi_{[0,+\infty)}(t) K Z(t)$. However, the resulting $u_D(t) = u(t-D) = \chi_{[D,+\infty)}(t) K Z(t-D)$ is discontinuous at $t=D$ while $u_D$ must be of class $\mathcal{C}^2$ over $\mathbb{R}_+$ to ensure the existence of a classical solution of (\ref{def: boundary control system}). Let $t_0 > 0$ be given. We consider a transition signal (from open loop to closed loop) $\varphi \in \mathcal{C}^2([-D,+\infty);\mathbb{R})$ which is such that $0 \leq \varphi \leq 1$, $\left. \varphi \right\vert_{[-D,0]} = 0$, and $\left. \varphi \right\vert_{[t_0,+\infty)} = 1$. We define the control input $u(t) = \varphi(t) K Z(t)$. It satisfies $\left. u \right\vert_{[-D,0]} = 0$ and, for all $t \geq 0$, 
\begin{align}
u(t) & = \varphi(t) K Z(t) \nonumber \\ 
	 & = \varphi(t) K Y(t) \label{eq: command input u} \\
	 & \phantom{=}\, + \varphi(t) K \int_{\max(t-D,0)}^{t} e^{(t-s-D)A_{N_0}} B_{N_0} u(s) \diff s , \nonumber
\end{align}
where it has been used that the system is uncontrolled for $t \leq 0$. 
In particular, the control law is such that $u_D(t) = u(t-D) = \varphi(t-D) K Z(t-D)$ with $u_D(t) = 0$ for $t \leq D$ and $u_D(t) = K Z(t-D)$ for $t \geq D+t_0$.

\subsection{Characterization of the control law}
In practice, it is convenient to use the control law expressed under the form (\ref{eq: command input u}) since it allows its computation at time $t$ based on the measure of $Y$ at time $t$ and the past history of the control law $u$. To do so, we must show that (\ref{eq: command input u}) fully characterizes $u$, i.e., the uniqueness of the function $u$ satisfying the implicit equation (\ref{eq: command input u}). In other words, it requires to invert the Artstein transformation~\cite{bresch2018new} when weighted by the transition signal $\varphi$. For any locally integrable function $f:\mathbb{R}_+ \rightarrow \mathbb{K}^m$, we define $T_{D}f:\mathbb{R}_+ \rightarrow \mathbb{K}^m$ as follows:
\begin{equation*}
(T_{D}f)(t) = \varphi(t) K \int_{\max(t-D,0)}^{t} e^{(t-s-D)A_{N_0}}B_{N_0}f(s) \diff s.
\end{equation*}
In particular $T_{D}f \in \mathcal{C}^0(\mathbb{R}_+;\mathbb{K}^m)$, and thus we can consider the iterations $T_{D}^{n}f$ for any $n \in \mathbb{N}$. 

\begin{lemma}\label{eq: inv weighted Artstein transformation}
Let $D>0$, $T \in \mathbb{R}_+^* \cup \{+\infty\}$, $g \in \mathcal{C}^0([0,T];\mathbb{K}^{N_0})$, and $\varphi \in \mathcal{C}^0([0,T];\mathbb{R})$ such that $0 \leq \varphi \leq 1$ be given. Then, there exists a unique locally integrable function $v$ defined over $[0,T]$ such that for all $t \in [0,T]$,
\begin{align*}
v(t) & = \varphi(t) K g(t) \\ 
& \phantom{=}\, + \varphi(t) K \int_{\max(t-D,0)}^{t} e^{(t-s-D)A_{N_0}} B_{N_0} v(s) \diff s . \nonumber
\end{align*}
Furthermore $v \in \mathcal{C}^0([0,T];\mathbb{K}^m)$ and is given by the following series expansion:
\begin{equation*}
v(t) = \sum\limits_{k=0}^{\infty} (T_{D}^{k}(\varphi K g))(t) ,
\end{equation*}
where the series converge uniformingly over any time interval of finite length.
\end{lemma}

The inversion of the Artstein transformation corresponding to the case $\varphi = 1$ has been investigated in~\cite{bresch2018new}. The proof of Lemma~\ref{eq: inv weighted Artstein transformation} is a straightfoward extension of the proof of Theorem~1 in~\cite{bresch2018new} by noting that $v = \varphi K g + T_D v$ and by using the fact that $\varphi$ is a continuous function that satisfies $0 \leq \varphi \leq 1$.

\section{Study of the closed-loop infinite-dimensional system}\label{sec: Study of the infinite-dimensional closed-loop system}
Throughout this section, we assume that Assumptions~\ref{assum: A1}, \ref{assum: A2}, and~\ref{assum: A3} hold true. Under these conditions, it has been proposed in Section~\ref{sec: construction of the feedback control strategy} to resort to the control law given by (\ref{eq: command input u}) to stabilize the infinite-dimensional system (\ref{def: boundary control system}). As the control law has been derived on a finite-dimensional part of the original infinite-dimensional system, we must guarantee that the proposed control strategy successfully stabilizes the full system. Furthermore, in order to make valid the computations performed in the previous section, we must ensure that the \emph{a priori} assumption is indeed satisfied, i.e., the proposed control law is of class $\mathcal{C}^2$.

\subsection{Dynamics of the closed-loop system}
Let $D,t_0 > 0$ be given. We consider a given transition signal $\varphi \in \mathcal{C}^2([-D,+\infty);\mathbb{R})$ such that $0 \leq \varphi \leq 1$, $\left. \varphi \right\vert_{[-D,0]} = 0$, and $\left. \varphi \right\vert_{[t_0,+\infty)} = 1$. The closed-loop system dynamics takes the following form:
\begin{equation}\label{def: boundary control system - closed-loop}
\left\{\begin{split}
\dfrac{\mathrm{d} X}{\mathrm{d} t}(t) & = \mathcal{A} X(t) + d(t) , \\
\mathcal{B} X (t) & = u_D(t) = u(t-D) , \\
\left. u \right\vert_{[-D,0]} & = 0 \\
u(t) & = \varphi(t) K Y(t) \\
& \phantom{=} + \varphi(t) K \int_{\max(t-D,0)}^{t} e^{(t-s-D)A_{N_0}} B_{N_0} u(s) \diff s , \\
Y(t) & = 
\begin{bmatrix}
\left< X(t) , \psi_1 \right>_\mathcal{H} \\ \vdots \\ \left< X(t) , \psi_{N_0} \right>_\mathcal{H}
\end{bmatrix} , \\
X(0) & = X_0 
\end{split}\right.
\end{equation}
for any $t \geq 0$. The feedback gain $K \in \mathbb{K}^{m \times N_0}$ is such that $A_\mathrm{cl} \triangleq A_{N_0} + e^{-D A_{N_0}} B_{N_0} K$ is Hurwitz (with desired pole placement). Function $d : \mathbb{R}_+ \rightarrow \mathcal{H}$ represents a distributed disturbance.

\subsection{Well-posedness in terms of classical solutions}
The following lemma ensures both the well-posedness of the closed-loop system in terms of classical solutions and the sufficient regularity of the control input.

\begin{lemma}\label{lemma: well-posedness closed-loop PDE}
Let $(\mathcal{A},\mathcal{B})$ be an abstract boundary control system such that Assumptions~\ref{assum: A1}, \ref{assum: A2}, and~\ref{assum: A3} hold true. For any $X_0 \in D(\mathcal{A}_0)$ and  $d \in \mathcal{C}^1(\mathbb{R}_+;\mathcal{H})$, the closed-loop system (\ref{def: boundary control system - closed-loop}) admits a unique classical solution $X \in \mathcal{C}^0(\mathbb{R}_+;D(\mathcal{A})) \cap \mathcal{C}^1(\mathbb{R}_+;\mathcal{H})$. The associated control law $u$ is uniquely defined and is of class $\mathcal{C}^2([-D,+\infty);\mathbb{K}^m)$. It can be written under the form $u = \varphi K Z$ with, for all $t \geq 0$,
\begin{equation}\label{eq: Artstein transform}
Z(t) \triangleq Y(t) + \int_{t-D}^{t} e^{(t-s-D)A_{N_0}} B_{N_0} u(s) \diff s ,
\end{equation}
which is such that $Z \in \mathcal{C}^2(\mathbb{R}_+;\mathbb{K}^{N_0})$ and satisfies, for all $t \geq 0$,
\begin{equation}\label{eq: EDO satisfied by Z}
\dot{Z}(t) = (A_{N_0} + \varphi(t) e^{-D A_{N_0}} B_{N_0} K) Z(t) + D_{N_0}(t) ,
\end{equation}
where $D_{N_0}(t)$ is defined by (\ref{eq: def DN0}). In particular, for all $t \geq D+t_0$,
\begin{equation}\label{eq: EDO satisfied by Z - t geq D+t0}
\dot{Z}(t) = A_\mathrm{cl} Z(t) + D_{N_0}(t) .
\end{equation}
Furthermore, $u$ is also expressed for all $t \geq 0$ by the following series expansion
\begin{equation}\label{eq: u expressed by series}
u(t) = \sum\limits_{k=0}^{\infty} (T_{D}^{k}(\varphi K Y))(t) ,
\end{equation}
where the series converges uniformingly over any time interval of finite length. 
\end{lemma}

\textbf{Proof.} 
We first note that, as $\left. u \right\vert_{[-D,0]} = 0$, (\ref{def: boundary control system - closed-loop}) is equivalent over the time interval $[0,D]$ to the following standard evolution problem 
\begin{equation*}
\left\{\begin{split}
\dfrac{\mathrm{d} X}{\mathrm{d} t}(t) & = \mathcal{A}_0 X(t) + d(t) , & t \in [0,D] \\
X(0) & = X_0 
\end{split}\right.
\end{equation*}
As $\mathcal{A}_0$ generates a $C_0$-semigroup, we deduce (see, e.g., \cite{Curtain2012}) the existence and the uniqueness of a classical solution $X \in \mathcal{C}^0([0,D];D(\mathcal{A})) \cap \mathcal{C}^1([0,D];\mathcal{H})$ such that (\ref{def: boundary control system - closed-loop}) holds true over the time interval $[0,D]$.

We now proceed by induction. Assume that, for a given $n \in \mathbb{N}^*$, there exists a unique classical solution over the time interval $[0,nD]$ denoted by $X \in \mathcal{C}^0([0,nD];D(\mathcal{A})) \cap \mathcal{C}^1([0,nD];\mathcal{H})$ of (\ref{def: boundary control system - closed-loop}) with associated control input $u \in \mathcal{C}^0([-D,(n-1)D])$ satisfying 
$\left. u \right\vert_{[-D,0]} = 0$ and, for all $0 \leq t \leq (n-1)D$,
\begin{equation}\label{eq: recursive - implit equation u}
u(t) = \varphi(t) K Y(t) + (T_{D}u)(t) .
\end{equation}
We show that there exists a unique classical solution $\tilde{X} \in \mathcal{C}^0([0,(n+1)D];D(\mathcal{A})) \cap \mathcal{C}^1([0,(n+1)D];\mathcal{H})$ of (\ref{def: boundary control system - closed-loop}) over the time interval $[0,(n+1)D]$ with a uniquely defined associated control input $\tilde{u}$. In particular, such a solution must satisfy (\ref{def: boundary control system - closed-loop}) over the restricted time interval $[0,nD]$ and thus, by induction hypothesis, we must have $\left. \tilde{X} \right\vert_{[0,nD]} = X$. Furthermore, $\tilde{X}$ must satisfy
\begin{equation}\label{eq: existence classical solutions tilde_X}
\left\{\begin{split}
\dfrac{\mathrm{d} \tilde{X}}{\mathrm{d} t}(t) & = \mathcal{A} \tilde{X}(t) + d(t) , & t \in [nD,(n+1)D] \\
\mathcal{B} \tilde{X} (t) & = \tilde{u}_D(t) = \tilde{u}(t-D) , & t \in [nD,(n+1)D] \\
\left. \tilde{u} \right\vert_{[-D,0]} & = 0 \\
\tilde{u}(t) & = \varphi(t) K Y(t) + (T_{D}\tilde{u})(t) , & t \in [0,nD] \\
Y(t) & = 
\begin{bmatrix}
\left< X(t) , \psi_1 \right>_\mathcal{H} \\ \vdots \\ \left< X(t) , \psi_{N_0} \right>_\mathcal{H}
\end{bmatrix} , & t \in [0,nD] \\
\tilde{X}(nD) & = X(nD) 
\end{split}\right.
\end{equation}
Note that, due to the delay $D > 0$, the control input $\tilde{u}$ is only defined by $X$ over the time interval $[0,nD]$ and does not depend on $\tilde{X}$ over $[nD,(n+1)D]$. As $X \in \mathcal{C}^0([0,nD];D(\mathcal{A})) \cap \mathcal{C}^1([0,nD];\mathcal{H})$, we  have that $Y \in \mathcal{C}^1([0,nD];\mathbb{K}^{N_0})$. Then, according to the Lemma~\ref{eq: inv weighted Artstein transformation}, 1) the control $\tilde{u}$ is well and uniquely defined over $[-D,nD]$; 2) $\tilde{u}$ is continuous over $[-D,nD]$; 3) as both $u$ and $\left. \tilde{u} \right\vert_{[-D,(n-1)D]}$ satisfy (\ref{eq: recursive - implit equation u}) for all $t \in [0,(n-1)D]$, we have by uniqueness that $\left. \tilde{u} \right\vert_{[-D,(n-1)D]} = u$. Furthermore, we can write $\tilde{u}(t) = \varphi(t) K Z(t)$ with, for all $t \in [0,nD]$,
\begin{equation*}
Z(t) = Y(t) + \int_{t-D}^{t} e^{(t-s-D)A_{N_0}} B_{N_0} \tilde{u}(s) \diff s .
\end{equation*}
Thus, we infer that $Z \in \mathcal{C}^1([0,nD];\mathbb{K}^{N_0})$. As $X$ is a classical solution of (\ref{def: boundary control system - closed-loop}) over the time interval $[0,nD]$, we obtain with the same approach used to derive (\ref{eq: ODE satisfies by Y}) that $Y$ satisfies the following ODE over the time interval $[0,nD]$
\begin{equation*}
\dot{Y}(t) 
= A_{N_0} Y(t) + B_{N_0} \tilde{u}(t-D) + D_{N_0}(t) ,
\end{equation*}
where $D_{N_0}(t)$ is defined by (\ref{eq: def DN0}). Thus, we have for all $t \in [0,nD]$,
\begin{align*}
\dot{Z}(t) & = A_{N_0} Z(t) + e^{-D A_{N_0}} B_{N_0} \tilde{u}(t) + D_{N_0}(t) \\
& = (A_{N_0} + \varphi(t) e^{-D A_{N_0}} B_{N_0} K) Z(t)  + D_{N_0}(t).
\end{align*}
As $d \in \mathcal{C}^1(\mathbb{R}_+;\mathcal{H})$, we have $D_{N_0} \in \mathcal{C}^1(\mathbb{R}_+;\mathbb{K}^{N_0})$. We deduce that $Z$ is of class $\mathcal{C}^2$ over $[0, nD]$. Thus, the control law satisfies $\tilde{u} = \varphi K Z \in \mathcal{C}^2([-D, nD];\mathbb{K}^{N_0})$, showing that $\tilde{u}_D \in \mathcal{C}^2([0, (n+1)D];\mathbb{K}^{N_0})$. Furthermore, the distributed disturbance is such that $d \in \mathcal{C}^1(\mathbb{R}_+;\mathcal{H})$ while the initial condition of (\ref{eq: existence classical solutions tilde_X}) given at $t = nD$ is such that $X(nD) \in D(\mathcal{A})$ and $\mathcal{B}X(nD) = u_D(nD) = \tilde{u}_D(nD)$. This yields (see, e.g., \cite[Th.~3.3.3]{Curtain2012}) the existence and uniqueness of the classical solution $\left. \tilde{X} \right\vert_{[nD,(n+1)D]}$ associated with (\ref{eq: existence classical solutions tilde_X}). As $\tilde{X}(nD) = X(nD) $ and $\dfrac{\mathrm{d} \tilde{X}}{\mathrm{d} t}(nD) = \mathcal{A} \tilde{X}(nD) = \mathcal{A} X(nD) = \dfrac{\mathrm{d} X}{\mathrm{d} t}(nD)$, it shows that the obtained $\tilde{X}$ is such that $\tilde{X}\in \mathcal{C}^0([0,(n+1)D];D(\mathcal{A})) \cap \mathcal{C}^1([0,(n+1)D];\mathcal{H})$ and is the unique classical solution of (\ref{def: boundary control system - closed-loop}) over $[0,(n+1)D]$.

By induction, it shows the existence and the uniqueness of a classical solutions $X \in \mathcal{C}^0(\mathbb{R}_+;D(\mathcal{A})) \cap \mathcal{C}^1(\mathbb{R}_+;\mathcal{H})$ for the closed-loop system (\ref{def: boundary control system - closed-loop}) associated with $X_0 \in D(\mathcal{A}_0)$ and $d \in \mathcal{C}^1(\mathbb{R}_+;\mathcal{H})$. The claimed properties of the control input $u$ directly follow from the above developments and the application of Lemma~\ref{eq: inv weighted Artstein transformation}. \qed

\subsection{Exponential ISS property of the closed-loop system}
This section is devoted to the demonstration of the following stability result.

\begin{theorem}\label{thm: ISS}
Let $(\mathcal{A},\mathcal{B})$ be an abstract boundary control system such that Assumptions~\ref{assum: A1}, \ref{assum: A2}, and~\ref{assum: A3} hold true. There exist constants $\overline{C}_1,\overline{C}_2,\overline{C}_3,\overline{C}_4 \in \mathbb{R}_+$ such that, for any $X_0 \in D(\mathcal{A}_0)$ and $d \in \mathcal{C}^1(\mathbb{R}_+;\mathcal{H})$, the classical solution solution $X$ of (\ref{def: boundary control system - closed-loop}) associated with the initial condition $X_0$ and the distributed disturbance $d$ satisfies the ISS estimate
\begin{equation}\label{eq: global exponential stability}
\Vert X(t) \Vert_\mathcal{H} 
\leq
\overline{C}_1 e^{- \kappa_0 t} \Vert X_0 \Vert_\mathcal{H}
+ \overline{C}_2 \sup\limits_{\tau \in [0,t]} \Vert d(\tau) \Vert_\mathcal{H} ,
\end{equation}
and the control law satisfies
\begin{equation}\label{eq: command effort global exponential convergence}
\Vert u(t) \Vert
\leq
\overline{C}_3 e^{- \kappa_0 t} \Vert X_0 \Vert_\mathcal{H}
+ \overline{C}_4 \sup\limits_{\tau \in [0,t]} \Vert d(\tau) \Vert_\mathcal{H} ,
\end{equation}
for all $t \geq 0$.
\end{theorem}

\begin{remark}
Theorem~\ref{thm: ISS} ensures the stability of the closed-loop system whatever the value of the delay $D > 0$ may be.
\end{remark}

To prove the theorem, we consider throughout this section $X_0 \in D(\mathcal{A}_0)$ and $d \in \mathcal{C}^1(\mathbb{R}_+;\mathcal{H})$ arbitrarily given. Let $X \in \mathcal{C}^0(\mathbb{R}_+;D(\mathcal{A})) \cap \mathcal{C}^1(\mathbb{R}_+;\mathcal{H})$ be the classical solution of the closed-loop system (\ref{def: boundary control system - closed-loop}) associated with the initial condition $X_0 \in D(\mathcal{A}_0)$ and the distributed disturbance $d \in \mathcal{C}^1(\mathbb{R}_+;\mathcal{H})$. We denote by $Z$ the function defined by (\ref{eq: Artstein transform}).

\subsubsection{Definition of the Lyapunov function candidate}
The proof of the theorem relies on the following Lyapunov function candidate, defined for $t \geq 0$ by
\begin{align}
V(t) 
& = \gamma_1 \left\{ Z(t)^* P Z(t) + \int_{t-D}^{t} \varphi(s) Z(s)^* P Z(s) \diff s \right\} \nonumber \\
& \phantom{=}\, + \gamma_2 \varphi(t-D) Z(t-D)^* P Z(t-D) \label{Lyapunov function} \\
& \phantom{=}\, + \dfrac{1}{2} \sum\limits_{k \geq N_0+1} \left\vert \left< X(t) - Bu_D(t) , \psi_k \right>_\mathcal{H} \right\vert^2 , \nonumber
\end{align}
where, because $A_\mathrm{cl} = A_{N_0} + e^{-D A_{N_0}} B_{N_0} K$ is Hurwitz, $P \in \mathcal{H}_{N_0}^{+*}$ is a Hermitian definite positive matrix such that 
\begin{equation}\label{eq: Lypunov identity}
A_\mathrm{cl}^* P + P A_\mathrm{cl} = -I_{N_0} .
\end{equation}
Constant $\gamma_1,\gamma_2 \in \mathbb{R}_+^*$ are sufficiently large parameters to be selected latter, independently of the initial condition $X_0$ and the distributed disturbance $d$. Note that, from the definition, one has $V(t) \geq 0$ for all $t \geq 0$. Thus the selection of $\gamma_1$ and $\gamma_2$ will be only driven to ensure the exponential decay of $V$.

\begin{remark}
Function $V$ is well-defined and belongs to $\mathcal{C}^1(\mathbb{R}_+;\mathbb{R})$. Indeed, as $\varphi$ and $Z$ are continuous over $\mathbb{R}_+$, the integral term is finite and, from (\ref{eq: Riesz basis - inequality}),
\begin{align*}
& \sum\limits_{k \geq N_0+1} \left\vert \left< X(t) - Bu_D(t) , \psi_k \right>_\mathcal{H} \right\vert^2 \\
& \qquad\leq 
\sum\limits_{k \geq 1} \left\vert \left< X(t) - Bu_D(t) , \psi_k \right>_\mathcal{H} \right\vert^2 \\
& \qquad\leq
\dfrac{1}{m_R} \Vert X(t) - Bu_D(t) \Vert_\mathcal{H}
< \infty .
\end{align*}
Thus we have $V(t) \in \mathbb{R}_+$. The continuous differentiability of $V$ follows from Annex~\ref{annex1} and the fact that functions $\varphi$, $X$, $Z$, and $u$ are of class $\mathcal{C}^1$.
\end{remark}

\begin{remark}
At this point, it is relevant to discuss the motivation behind the choice of the different terms of the Lyapunov function candidate (\ref{Lyapunov function}). 
\begin{enumerate}
\item Assuming a zero distributed disturbance ($d = 0$), the term $Z(t)^* P Z(t)$ provides, based on (\ref{eq: Lypunov identity}), a Lyapunov function for the finite-dimensional system $\dot{Z}(t) = A_\mathrm{cl} Z(t)$. It aims at ensuring the exponential convergence to zero of the $N_0$ first coefficients $\left< X(t) , \psi_n \right>_\mathcal{H}$ corresponding to the projection of the system trajectory $X(t)$ into the Riesz basis $\left\{ \phi_n , \; n \in \mathbb{N}^* \right\}$ (see (\ref{eq: proj Riesz basis})). 
\item In order to ensure the stability of the full infinite-dimensional system, the Lyapunov function candidate $V$ must ensure the convergence of all the modes, including the coefficients $\left< X(t) , \psi_n \right>_\mathcal{H}$, $n \geq N_0 +1$, which were not considered in the synthesis of the control law. A natural choice to capture these coefficients would consist in the use of the term $\dfrac{1}{2} \sum\limits_{k \geq N_0+1} \left\vert \left< X(t) , \psi_k \right>_\mathcal{H} \right\vert^2$. However, the ODE describing the time domain evolution of $\left< X(t) , \psi_n \right>_\mathcal{H}$ given by (\ref{eq: coeff in Riesz basis ODE}) shows that the eigenvalue $\lambda_n$ appears via the following term: $\lambda_n \left< X(t) - B u_D(t) , \psi_n \right>_\mathcal{H}$. Therefore, in order to be able to absorb all the occurrences of the eigenvalue $\lambda_n$, $n \geq N_0 + 1$, via the inequality $\operatorname{Re} \lambda_n \leq - \alpha$ of Assumption~\ref{assum: A2}, we consider the term $\dfrac{1}{2} \sum\limits_{k \geq N_0+1} \left\vert \left< X(t) - Bu_D(t) , \psi_k \right>_\mathcal{H} \right\vert^2$ (see (\ref{eq: full absorbtion lambda_k}) for details).
\item As $u = \varphi K Z$, the introduction of the term $u_D(t)$ in the Lyapunov function candidate $V$ yields the occurrence of the term $Z(t-D)$. It requires the introduction of the term $\varphi(t-D) Z(t-D)^* P Z(t-D)$ for compensation purposes. The switching signal $\varphi$ is used to materialize the fact that the contribution of this term is relevant only for $t \geq D$.
\item Finally, the contribution of the term $\int_{t-D}^{t} \varphi(s) Z(s)^* P Z(s) \diff s$ is to provide an upper bound on the norm of the system trajectory $X(t)$ which only depends on $V(t)$ (see Lemma~\ref{eq: estimate X by above with V}). 
\end{enumerate} 
\end{remark}
The detailed properties of the Lyapunov function candidate $V$ are detailed in the next lemmas.

\subsubsection{Upper bound on the norm of $X$}

First, we establish a connection between the norm of the system trajectory $X(t)$ and the value of the Lyapunov function candidate $V(t)$. We define the constant $C_1 > 0$ by
\begin{equation}\label{eq: constant C1}
C_1 \triangleq 2 \max \left( 1 , D e^{2 D \left\Vert A_{N_0} \right\Vert} \left\Vert B_{N_0} K \right\Vert^2 \right) .
\end{equation}
We denote by $\lambda_m(P) > 0$ the smallest eigenvalue of $P$.

\begin{lemma}\label{eq: estimate X by above with V}
Under the assumptions of Theorem~\ref{thm: ISS} and for $\gamma_1 > C_1/\lambda_m(P)$ and $\gamma_2 > \Vert B K \Vert^2/(m_R \lambda_m(P))$ arbitrarily given, there exists a constant $C_4 = C_4(\gamma_2) > 0$, independent of $X_0$ and $d$, such that 
\begin{equation}\label{eq: X bounded above by cste times sqrt V}
\left\Vert X(t) \right\Vert_\mathcal{H} \leq C_4 \sqrt{V(t)} 
\end{equation}
for all $t \geq 0$.
\end{lemma}

\textbf{Proof.} 
From (\ref{eq: Artstein transform}) and using the identity $u = \varphi K Z$, we have that for all $t \geq 0$,
\begin{equation*}
Y(t) = Z(t) - \int_{t-D}^{t} \varphi(s) e^{(t-s-D)A_{N_0}} B_{N_0} K Z(s) \diff s .
\end{equation*}
Using the Cauchy-Schwartz (C.S.) inequality and the fact that $0 \leq \varphi \leq 1$, we deduce that, for all $t \geq 0$,
\begin{align*}
& \left\Vert Y(t) \right\Vert \\
& \leq \left\Vert Z(t) \right\Vert + \left\Vert \int_{t-D}^{t} \varphi(s) e^{(t-s-D)A_{N_0}} B_{N_0} K Z(s) \diff s \right\Vert \\
     & \leq \left\Vert Z(t) \right\Vert + e^{D \left\Vert A_{N_0} \right\Vert} \left\Vert B_{N_0} K \right\Vert \int_{t-D}^{t} \varphi(s) \left\Vert Z(s) \right\Vert \diff s \\
     & \overset{C.S.}{\leq} \left\Vert Z(t) \right\Vert + \sqrt{D} e^{D \left\Vert A_{N_0} \right\Vert} \left\Vert B_{N_0} K \right\Vert \sqrt{\int_{t-D}^{t} \varphi(s) \left\Vert Z(s) \right\Vert^2 \diff s} ,
\end{align*}
which gives 
\begin{align}
& \left\Vert Y(t) \right\Vert^2 \nonumber \\
& \leq 2 \left\Vert Z(t) \right\Vert^2 
+  2 D e^{2 D \left\Vert A_{N_0} \right\Vert} \left\Vert B_{N_0} K \right\Vert^2 \int_{t-D}^{t} \varphi(s) \left\Vert Z(s) \right\Vert^2 \diff s \nonumber \\
& \leq C_1 \left\{ \left\Vert Z(t) \right\Vert^2 + \int_{t-D}^{t} \varphi(s) \left\Vert Z(s) \right\Vert^2 \diff s \right\} , \label{eq: norm Y(t) upper bound}
\end{align}
where $C_1$ is defined by (\ref{eq: constant C1}). Now, from the definition of $V$ given by (\ref{Lyapunov function}) and using (\ref{eq: Riesz basis - inequality}), we have for all $t \geq 0$,
\begin{align*}
V(t) & 
\geq \gamma_1 \lambda_m(P) \left\{ \left\Vert Z(t) \right\Vert^2 + \int_{t-D}^{t} \varphi(s) \left\Vert Z(s) \right\Vert^2 \diff s \right\} \\
& \phantom{\geq}\; + \gamma_2 \lambda_m(P) \varphi(t-D) \left\Vert Z(t-D) \right\Vert^2 \\
& \phantom{\geq}\; + \dfrac{1}{2 M_R} \left\Vert X(t) - Bu_D(t) \right\Vert_\mathcal{H}^2  \\
& \phantom{\geq}\; - \dfrac{1}{2} \sum\limits_{k = 1}^{N_0} \left\vert \left< X(t) - Bu_D(t) , \psi_k \right>_\mathcal{H} \right\vert^2 .
\end{align*}
Recalling that $u_D(t) = u(t-D) = \varphi(t-D) K Z(t-D)$ and $0 \leq \varphi \leq 1$ which gives $\varphi^2 \leq \varphi$, we have
\begin{align*}
& \sum\limits_{k = 1}^{N_0} \left\vert \left< X(t) - Bu_D(t) , \psi_k \right>_\mathcal{H} \right\vert^2 \\
& \qquad\leq 2 \sum\limits_{k = 1}^{N_0} \left\{ \left\vert \left< X(t) , \psi_k \right>_\mathcal{H} \right\vert^2 + \left\vert \left< Bu_D(t) , \psi_k \right>_\mathcal{H} \right\vert^2 \right\} \\
& \qquad\leq 2  \Vert Y(t) \Vert^2 + 2 \sum\limits_{k \geq 1} \left\vert \left< Bu_D(t) , \psi_k \right>_\mathcal{H} \right\vert^2 \\
& \qquad\overset{(\ref{eq: Riesz basis - inequality})}{\leq} 2  \Vert Y(t) \Vert^2 + \dfrac{2}{m_R} \Vert Bu_D(t)\Vert_\mathcal{H}^2 \\
& \qquad\leq 2  \Vert Y(t) \Vert^2 + \dfrac{2 \Vert B K \Vert^2}{m_R} \{\varphi(t-D)\}^2 \Vert Z(t-D) \Vert^2 \\
& \qquad\leq 2  \Vert Y(t) \Vert^2 + \dfrac{2 \Vert B K \Vert^2}{m_R} \varphi(t-D) \Vert Z(t-D) \Vert^2 .
\end{align*}
We deduce that 
\begin{align*}
V(t) & 
\geq \gamma_1 \lambda_m(P) \left\{ \left\Vert Z(t) \right\Vert^2 + \int_{t-D}^{t} \varphi(s) \left\Vert Z(s) \right\Vert^2 \diff s \right\} - \Vert Y(t) \Vert^2 \\
& \phantom{\geq}\; + \left\{ \gamma_2 \lambda_m(P) - \dfrac{\Vert B K \Vert^2}{m_R} \right\} \varphi(t-D) \left\Vert Z(t-D) \right\Vert^2 \\
& \phantom{\geq}\; + \dfrac{1}{2 M_R} \left\Vert X(t) - Bu_D(t) \right\Vert_\mathcal{H}^2 .
\end{align*}
Using (\ref{eq: norm Y(t) upper bound}), this yields for all $t \geq 0$,
\begin{align*}
V(t) & 
\geq \left\{ \gamma_1 \lambda_m(P) - C_1 \right\} \left\{ \left\Vert Z(t) \right\Vert^2 + \int_{t-D}^{t} \varphi(s) \left\Vert Z(s) \right\Vert^2 \diff s \right\} \\
& \phantom{\geq}\; + \left\{ \gamma_2 \lambda_m(P) - \dfrac{\Vert B K \Vert^2}{m_R} \right\} \varphi(t-D) \left\Vert Z(t-D) \right\Vert^2 \\
& \phantom{\geq}\; + \dfrac{1}{2 M_R} \left\Vert X(t) - Bu_D(t) \right\Vert_\mathcal{H}^2 .
\end{align*}
As $\gamma_1,\gamma_2 \in \mathbb{R}_+^*$ are such that $\gamma_1 > C_1/\lambda_m(P)$ and $\gamma_2 > \Vert B K \Vert^2/(m_R \lambda_m(P))$, we have $C_2(\gamma_1) \triangleq \gamma_1 \lambda_m(P) - C_1 > 0$ and $C_3(\gamma_2) \triangleq \gamma_2 \lambda_m(P) - \dfrac{\Vert B K \Vert^2}{m_R} > 0$ are such that, for all $t \geq 0$,
\begin{align}
V(t) & 
\geq C_2(\gamma_1) \left\{ \left\Vert Z(t) \right\Vert^2 + \int_{t-D}^{t} \varphi(s) \left\Vert Z(s) \right\Vert^2 \diff s \right\} \nonumber \\
& \phantom{\geq}\, + C_3(\gamma_2) \varphi(t-D) \left\Vert Z(t-D) \right\Vert^2 \label{eq: lower bound V} \\ 
& \phantom{\geq}\, + \dfrac{1}{2 M_R} \left\Vert X(t) - Bu_D(t) \right\Vert_\mathcal{H}^2 . \nonumber
\end{align}
In particular, this yields for all $t \geq 0$,
\begin{align*}
\left\Vert X(t) \right\Vert_\mathcal{H} 
& \leq \left\Vert X(t) - Bu_D(t) \right\Vert_\mathcal{H} + \left\Vert Bu_D(t) \right\Vert_\mathcal{H}  \\
& \leq \sqrt{2 M_R V(t)} + \Vert BK \Vert \times \varphi(t-D) \left\Vert Z(t-D) \right\Vert \\
& \leq \sqrt{2 M_R V(t)} + \Vert BK \Vert \times \sqrt{\varphi(t-D)} \left\Vert Z(t-D) \right\Vert \\
& \leq \sqrt{2 M_R V(t)} + \Vert BK \Vert \times \dfrac{1}{\sqrt{C_3(\gamma_2)}} \sqrt{V(t)} \\
& \leq \left\{ \sqrt{2 M_R} + \dfrac{\Vert BK \Vert}{\sqrt{C_3(\gamma_2)}} \right\} \sqrt{V(t)} .
\end{align*}
Introducing $C_4 \triangleq \sqrt{2 M_R} + \dfrac{\Vert BK \Vert}{\sqrt{C_3(\gamma_2)}} > 0$, the claimed inequality (\ref{eq: X bounded above by cste times sqrt V}) holds true. \qed

\subsubsection{Exponential convergence of the closed-loop system trajectories}
In order to study the exponential decay of $V$, we consider the time interval over which the infinite-dimensional system is fully placed in closed loop, i.e., for $t > D+t_0$ which corresponds to $\varphi(t) = 1$. For $t > D+t_0$, one has
\begin{align*}
V(t) 
& = \gamma_1 \left\{ Z(t)^* P Z(t) + \int_{t-D}^{t} Z(s)^* P Z(s) \diff s \right\} \\
& \phantom{=}\, + \gamma_2 Z(t-D)^* P Z(t-D) \\
& \phantom{=}\, + \dfrac{1}{2} \sum\limits_{k \geq N_0+1} \left\vert \left< X(t) - Bu_D(t) , \psi_k \right> \right\vert^2 
\end{align*}
with $u_D(t) = u(t-D) = K Z(t-D)$. We also introduce the positive constant
\begin{equation}\label{eq: constant C5}
C_5 \triangleq \dfrac{2m}{\alpha m_R} \sum\limits_{i=1}^{m} \left\{ \Vert \mathcal{A} B e_i \Vert_\mathcal{H}^2 \Vert K_i \Vert^2 + \Vert B e_i \Vert_\mathcal{H}^2 \Vert K_i A_\mathrm{cl} \Vert^2 \right\} ,
\end{equation}
where $K_i$ is the $i$-th line of the  matrix of feedback gain $K$.

\begin{lemma}
Let $\beta \in (0,1)$ be arbitrarily given. Under the assumptions of Theorem~\ref{thm: ISS}, and for any arbitrarily given $\gamma_1 > C_1/\lambda_m(P)$ and $\gamma_2 > \max \left( \Vert B K \Vert^2/(m_R \lambda_m(P)) , C_5/(1-\beta) \right)$, there exist constants $\kappa_0 = \kappa_0(\beta,\gamma_2)>0$ and $C_6=(\beta,\gamma_1,\gamma_2) > 0$, independent of $X_0$ and $d$, such that we have for all $t \geq D+t_0$,
\begin{align}
\left\Vert X(t) \right\Vert_\mathcal{H} & \leq C_4 e^{- \kappa_0 (t-D-t_0)} \sqrt{V(D+t_0)} \label{eq: trajectories converge exponentially to zero} \\
& \phantom{\leq}\, + C_4 \sqrt{\dfrac{C_6}{2 \kappa_0}} \underset{\tau \in [0,t]}{\mathrm{sup}} \Vert d(\tau) \Vert_\mathcal{H} \nonumber
\end{align}
with a control input such that 
\begin{align}
\Vert u(t) \Vert
& \leq \dfrac{\Vert K \Vert}{\sqrt{C_2(\gamma_1)}} e^{- \kappa_0 (t-D-t_0)} \sqrt{V(D+t_0)} \nonumber \\
& \phantom{\leq}\, + \dfrac{\Vert K \Vert}{\sqrt{C_2(\gamma_1)}} \sqrt{\dfrac{C_6}{2 \kappa_0}} \underset{\tau \in [0,t]}{\mathrm{sup}} \Vert d(\tau) \Vert_\mathcal{H} . \label{eq: command effort converge exponentially to zero}
\end{align}
\end{lemma}

\textbf{Proof.} 
From the definition of $P$, we have that for all $t > t_0$,
\begin{align*}
\dfrac{\mathrm{d}}{\mathrm{d}t} \left[ Z^* P Z \right](t) 
& \overset{(\ref{eq: EDO satisfied by Z - t geq D+t0})}{=} Z(t)^* \left[ A_\mathrm{cl}^* P + P A_\mathrm{cl} \right] Z(t) \\
& \phantom{=}\, + D_{N_0}(t)^* P Z(t) + Z(t)^* P D_{N_0}(t) \\
& \overset{(\ref{eq: Lypunov identity})}{=} - \Vert Z(t) \Vert^2 + D_{N_0}(t)^* P Z(t) + Z(t)^* P D_{N_0}(t) .
\end{align*}
Thus, for all $t > D+t_0$,
\begin{align*}
& \dfrac{\mathrm{d}}{\mathrm{d}t} \left[ \int_{t-D}^{t} Z(s)^* P Z(s) \diff s \right] (t) \\
& \qquad =  Z(t)^* P Z(t) - Z(t-D)^* P Z(t-D) \\
& \qquad = - \int_{t-D}^{t} \Vert Z(s) \Vert^2 \diff s \\
& \qquad \phantom{=}\, + \int_{t-D}^{t} D_{N_0}(s)^* P Z(s) + Z(s)^* P D_{N_0}(s) \diff s.
\end{align*}
Let $\beta \in (0,1)$ be arbitrarily given. We infer from the Young inequality (Y.I.) that, for all $t > t_0$,
\begin{align*}
& \dfrac{\mathrm{d}}{\mathrm{d}t} \left[ Z^* P Z \right](t) \\
& \quad\leq - \Vert Z(t) \Vert^2 + 2 \Vert P \Vert \Vert D_{N_0}(t) \Vert \Vert Z(t) \Vert \\
& \quad\overset{\mathrm{Y.I.}}{\leq} - \Vert Z(t) \Vert^2 + 2 \left( \dfrac{\beta}{2} \Vert Z(t) \Vert^2 + \dfrac{1}{2 \beta} \Vert P \Vert^2 \Vert D_{N_0}(t) \Vert^2 \right) \\ 
& \quad\overset{(\ref{eq: above estimate of D_N0})}{\leq} - (1 - \beta) \Vert Z(t) \Vert^2 + \dfrac{\Vert P \Vert^2}{\beta m_R} \Vert d(t) \Vert_\mathcal{H}^2 ,
\end{align*}
and for all $t > D+t_0$,
\begin{align*}
& \dfrac{\mathrm{d}}{\mathrm{d}t} \left[ \int_{t-D}^{t} Z(s)^* P Z(s) \diff s \right](t) \\
& \leq - \int_{t-D}^{t} \Vert Z(s) \Vert^2 \diff s 
+ 2 \int_{t-D}^{t} \Vert P \Vert \Vert D_{N_0}(s) \Vert \Vert Z(s) \Vert \diff s \\
& \overset{\mathrm{Y.I.}}{\leq} - \int_{t-D}^{t} \Vert Z(s) \Vert^2 \diff s \\
& \phantom{\overset{\mathrm{Y.I.}}{\leq}}\, + 2 \int_{t-D}^{t} \dfrac{\beta}{2} \Vert Z(s) \Vert^2 + \dfrac{1}{2 \beta} \Vert P \Vert^2 \Vert D_{N_0}(s) \Vert^2 \diff s \\
& \overset{(\ref{eq: above estimate of D_N0})}{\leq} - (1-\beta) \int_{t-D}^{t} \Vert Z(s) \Vert^2 \diff s 
+ \dfrac{\Vert P \Vert^2}{\beta m_R} \int_{t-D}^{t} \Vert d(s) \Vert_\mathcal{H}^2 \diff s \\
& \leq - (1-\beta) \int_{t-D}^{t} \Vert Z(s) \Vert^2 \diff s 
+ \dfrac{D \Vert P \Vert^2}{\beta m_R} \underset{\tau \in [t-D,t]}{\mathrm{sup}} \Vert d(\tau) \Vert_\mathcal{H}^2 .
\end{align*}
Finally, we have (see Annex~\ref{annex1})
\begin{align*}
& \dfrac{\mathrm{d}}{\mathrm{d}t} \left[ \dfrac{1}{2} \sum\limits_{k \geq N_0+1} \left\vert \left< X - Bu_D , \psi_k \right>_\mathcal{H} \right\vert^2 \right](t) \\
& \qquad = \sum\limits_{k \geq N_0+1} \operatorname{Re} \left\{ 
\left< \dfrac{\mathrm{d}X}{\mathrm{d}t}(t) - B\dot{u}_D(t) , \psi_k \right>_\mathcal{H} \right. \\
& \qquad \phantom{= \sum\limits_{k \geq N_0+1} \mathrm{Re} \{}\;\; \times \left.
\overline{\left< X(t) - Bu_D(t) , \psi_k \right>_\mathcal{H}} 
\right\} .
\end{align*}
As $X$ is a classical solution of the abstract Cauchy problem, using (\ref{eq: coeff in Riesz basis ODE}), Assumption~\ref{assum: A1}, and the Young inequality, we have for $k \geq N_0+1$ that
\begin{align}
& \operatorname{Re} \left\{ 
\left< \dfrac{\mathrm{d}X}{\mathrm{d}t}(t) - B\dot{u}_D(t) , \psi_k \right>_\mathcal{H}
\overline{\left< X(t) - Bu_D(t) , \psi_k \right>_\mathcal{H}} 
\right\} \nonumber \\
& \overset{(\ref{eq: coeff in Riesz basis ODE})}{=}
\operatorname{Re} (\lambda_k) \left\vert \left< X(t) - B u_D(t) , \psi_k \right>_\mathcal{H} \right\vert^2 \label{eq: full absorbtion lambda_k} \\
& \phantom{=}\, + \operatorname{Re} \left\{ 
\left( \left< \mathcal{A} B u_D(t) , \psi_k \right>_\mathcal{H}
+ \left< d(t) , \psi_k \right>_\mathcal{H}
- \left< B\dot{u}_D(t) , \psi_k \right>_\mathcal{H} \right) \right. \nonumber \\
& \phantom{= \, + \operatorname{Re} \{ }\, \times \left. \overline{\left< X(t) - Bu_D(t) , \psi_k \right>_\mathcal{H}} 
\right\} \nonumber \\
& \leq - \alpha \left\vert \left< X(t) - B u_D(t) , \psi_k \right>_\mathcal{H} \right\vert^2 \nonumber \\
& \phantom{=}\, + \left\{ \left\vert \left< \mathcal{A} B u_D(t) , \psi_k \right>_\mathcal{H} \right\vert
+ \left\vert d_k(t) \right\vert
+ \left\vert \left< B\dot{u}_D(t) , \psi_k \right>_\mathcal{H} \right\vert \right\} \nonumber \\
& \phantom{= \, +}\, \times \left\vert \left< X(t) - Bu_D(t) , \psi_k \right>_\mathcal{H} \right\vert \nonumber \\
& \overset{\mathrm{Y.I.}}{\leq} - \dfrac{\alpha}{2} \left\vert \left< X(t) - B u_D(t) , \psi_k \right>_\mathcal{H} \right\vert^2 \nonumber \\
& \phantom{=}\, + \dfrac{1}{2\alpha} \left\{ \left\vert \left< \mathcal{A} B u_D(t) , \psi_k \right>_\mathcal{H} \right\vert
+ \left\vert d_k(t) \right\vert
+ \left\vert \left< B\dot{u}_D(t) , \psi_k \right>_\mathcal{H} \right\vert \right\}^2 \nonumber .
\end{align}
Introducing $K_i$ the $i$-th line of the  matrix of feedback gain $K$, one has, for all $t > D+t_0$,
\begin{equation*}
u_D(t) = u(t-D) = K Z(t-D)
= \sum\limits_{i=1}^{m} \left\{ K_i Z(t-D) \right\} e_i 
\end{equation*}
and
\begin{align*}
\dot{u}_D(t) 
& = \dot{u}(t-D) 
= K \dot{Z}(t-D) \\
& \overset{(\ref{eq: EDO satisfied by Z - t geq D+t0})}{=} K ( A_\mathrm{cl} Z(t-D) + D_{N_0}(t-D) ) \\
& = \sum\limits_{i=1}^{m} \left\{ K_i A_\mathrm{cl} Z(t-D) \right\} e_i + K D_{N_0}(t-D) .
\end{align*}
This yields
\begin{align*}
& \operatorname{Re} \left\{ 
\left< \dfrac{\mathrm{d}X}{\mathrm{d}t}(t) - B\dot{u}_D(t) , \psi_k \right>_\mathcal{H}
\overline{\left< X(t) - Bu_D(t) , \psi_k \right>_\mathcal{H}} 
\right\} \\
& \leq - \dfrac{\alpha}{2} \left\vert \left< X(t) - B u_D(t) , \psi_k \right>_\mathcal{H} \right\vert^2 \\
& \phantom{\leq}\, + \dfrac{1}{2\alpha} \left\{ \left\vert \left< \mathcal{A} B K Z(t-D) , \psi_k \right>_\mathcal{H} \right\vert 
+ \left\vert \left< B K A_\mathrm{cl} Z(t-D) , \psi_k \right>_\mathcal{H} \right\vert \right. \\
& \phantom{\leq}\, \phantom{+ \dfrac{1}{2\alpha} \{}\; \left. + \left\vert d_k(t) \right\vert
+ \left\vert \left< B K D_{N_0}(t-D) , \psi_k \right>_\mathcal{H} \right\vert \right\}^2 \\
& \leq - \dfrac{\alpha}{2} \left\vert \left< X(t) - B u_D(t) , \psi_k \right>_\mathcal{H} \right\vert^2 \\
& \phantom{\leq}\, + \dfrac{2}{\alpha} \left\{ \left\vert \left< \mathcal{A} B K Z(t-D) , \psi_k \right>_\mathcal{H} \right\vert^2
+ \left\vert \left< B K A_\mathrm{cl} Z(t-D) , \psi_k \right>_\mathcal{H} \right\vert^2 \right. \\
& \phantom{\leq}\, \phantom{+ \dfrac{2}{\alpha} \{}\; \left. + \left\vert d_k(t) \right\vert^2
+ \left\vert \left< B K D_{N_0}(t-D) , \psi_k \right>_\mathcal{H} \right\vert^2 \right\} \\
& \leq - \dfrac{\alpha}{2} \left\vert \left< X(t) - B u_D(t) , \psi_k \right>_\mathcal{H} \right\vert^2 \\
& \phantom{\leq}\, + \dfrac{2}{\alpha} \left\vert \sum\limits_{i=1}^{m} \left< \mathcal{A} B e_i , \psi_k \right>_\mathcal{H} K_i Z(t-D) \right\vert^2 \\
& \phantom{\leq}\, + \dfrac{2}{\alpha} \left\vert \sum\limits_{i=1}^{m} \left< B e_i , \psi_k \right>_\mathcal{H} K_i A_\mathrm{cl} Z(t-D) \right\vert^2 \\
& \phantom{\leq}\, + \dfrac{2}{\alpha} \left\{ 
\left\vert d_k(t) \right\vert^2
+ \left\vert \left< B K D_{N_0}(t-D) , \psi_k \right>_\mathcal{H} \right\vert^2 \right\} \\
& \leq - \dfrac{\alpha}{2} \left\vert \left< X(t) - B u_D(t) , \psi_k \right>_\mathcal{H} \right\vert^2 \\
& \phantom{\leq}\, + \dfrac{2m}{\alpha} \sum\limits_{i=1}^{m} \left\vert \left< \mathcal{A} B e_i , \psi_k \right>_\mathcal{H} K_i Z(t-D) \right\vert^2 \\
& \phantom{\leq}\, + \dfrac{2m}{\alpha} \sum\limits_{i=1}^{m} \left\vert \left< B e_i , \psi_k \right>_\mathcal{H} K_i A_\mathrm{cl} Z(t-D) \right\vert^2 \\
& \phantom{\leq}\, + \dfrac{2}{\alpha} \left\{ 
\left\vert d_k(t) \right\vert^2
+ \left\vert \left< B K D_{N_0}(t-D) , \psi_k \right>_\mathcal{H} \right\vert^2 \right\} \\
& \leq - \dfrac{\alpha}{2} \left\vert \left< X(t) - B u_D(t) , \psi_k \right>_\mathcal{H} \right\vert^2 \\
& \phantom{\leq}\, + \dfrac{2m}{\alpha} \left\{ \sum\limits_{i=1}^{m} \left\vert \left< \mathcal{A} B e_i , \psi_k \right>_\mathcal{H}\right\vert^2 \Vert K_i \Vert^2 \right. \\
& \phantom{\leq}\, \phantom{+ \dfrac{2m}{\alpha} \{}\; \left. + \sum\limits_{i=1}^{m} \left\vert \left< B e_i , \psi_k \right>_\mathcal{H}\right\vert^2 \Vert K_i A_\mathrm{cl} \Vert^2 \right\} 
\Vert Z(t-D) \Vert^2 \\
& \phantom{\leq}\, + \dfrac{2}{\alpha} \left\{ 
\left\vert d_k(t) \right\vert^2
+ \left\vert \left< B K D_{N_0}(t-D) , \psi_k \right>_\mathcal{H} \right\vert^2 \right\} .
\end{align*}
We deduce that, for $t > D+t_0$,
\begin{align*}
& \dfrac{\mathrm{d}}{\mathrm{d}t} \left[ \dfrac{1}{2} \sum\limits_{k \geq N_0+1} \left\vert \left< X - Bu_D , \psi_k \right>_\mathcal{H} \right\vert^2 \right](t) \\
& \leq - \dfrac{\alpha}{2} \sum\limits_{k \geq N_0+1} \left\vert \left< X(t) - B u_D(t) , \psi_k \right>_\mathcal{H} \right\vert^2 \\
& \phantom{\leq}\, + \dfrac{2m}{\alpha} \sum\limits_{k \geq N_0+1} \left\{ \sum\limits_{i=1}^{m} \left\vert \left< \mathcal{A} B e_i , \psi_k \right>_\mathcal{H}\right\vert^2 \Vert K_i \Vert^2 \right. \\
& \phantom{\leq}\, \phantom{\dfrac{2m}{\alpha} \sum\limits_{k \geq N_0+1} \{} \left. + \sum\limits_{i=1}^{m} \left\vert \left< B e_i , \psi_k \right>_\mathcal{H}\right\vert^2 \Vert K_i A_\mathrm{cl} \Vert^2 \right\} 
\Vert Z(t-D) \Vert^2 \\
& \phantom{\leq}\, + \dfrac{2}{\alpha} \sum\limits_{k \geq N_0+1} \left\{ 
\left\vert d_k(t) \right\vert^2
+ \left\vert \left< B K D_{N_0}(t-D) , \psi_k \right>_\mathcal{H} \right\vert^2 \right\} \\
& \leq - \dfrac{\alpha}{2} \sum\limits_{k \geq N_0+1} \left\vert \left< X(t) - B u_D(t) , \psi_k \right>_\mathcal{H} \right\vert^2 \\
& \phantom{\leq}\, + \dfrac{2m}{\alpha} \left\{ \sum\limits_{i=1}^{m} \sum\limits_{k \geq 1} \left\vert \left< \mathcal{A} B e_i , \psi_k \right>_\mathcal{H}\right\vert^2 \Vert K_i \Vert^2 \right. \\
& \phantom{\leq}\, \left. \phantom{+ \dfrac{2m}{\alpha} \{}\;\; + \sum\limits_{i=1}^{m} \sum\limits_{k \geq 1} \left\vert \left< B e_i , \psi_k \right>_\mathcal{H}\right\vert^2 \Vert K_i A_\mathrm{cl} \Vert^2 \right\}
\Vert Z(t-D) \Vert^2 \\
& \phantom{\leq}\, + \dfrac{2}{\alpha} \sum\limits_{k \geq 1} \left\vert d_k(t) \right\vert^2
+ \dfrac{2}{\alpha} \sum\limits_{k \geq 1} \left\vert \left< B K D_{N_0}(t-D) , \psi_k \right>_\mathcal{H} \right\vert^2 \\
& \overset{(\ref{eq: Riesz basis - inequality})}{\leq} - \dfrac{\alpha}{2} \sum\limits_{k \geq N_0+1} \left\vert \left< X(t) - B u_D(t) , \psi_k \right>_\mathcal{H} \right\vert^2 
+ C_5 \Vert Z(t-D) \Vert^2 \\
& \phantom{\leq}\, + \dfrac{2}{\alpha m_R} \Vert d(t) \Vert_\mathcal{H}^2
+ \dfrac{2 \Vert BK \Vert^2}{\alpha m_R^2} \Vert d(t-D) \Vert_\mathcal{H}^2
\end{align*}
with constant $C_5$ given by (\ref{eq: constant C5}). As $\gamma_2 > C_5/(1-\beta)$, we deduce that, for all $t > D+t_0$,
\begin{align*}
\dot{V}(t) \leq & 
- \gamma_1 (1-\beta)
\left\{ \Vert Z(t) \Vert^2 + \int_{t-D}^{t} \Vert Z(s) \Vert^2 \diff s \right\} \\
& - (\gamma_2 (1-\beta)  - C_5) \Vert Z(t-D) \Vert^2 \\
& - \dfrac{\alpha}{2} \sum\limits_{k \geq N_0+1} \left\vert \left< X(t) - B u_D(t) , \psi_k \right>_\mathcal{H} \right\vert^2 \\
& + \dfrac{1}{m_R} \left( \dfrac{2}{\alpha} + \dfrac{\gamma_1 \Vert P \Vert^2}{\beta} \right) \Vert d(t) \Vert_\mathcal{H}^2 \\
& + \dfrac{1}{m_R} \left( \dfrac{2 \Vert BK \Vert^2}{\alpha m_R} + \dfrac{\gamma_2 \Vert P \Vert^2}{\beta} \right) \Vert d(t-D) \Vert_\mathcal{H}^2 \\ 
& + \dfrac{\gamma_1 D \Vert P \Vert^2}{\beta m_R} \underset{\tau \in [t-D,t]}{\mathrm{sup}} \Vert d(\tau) \Vert_\mathcal{H}^2 \\
\leq & 
- \dfrac{\gamma_1 (1-\beta)}{\lambda_M(P)}
\left\{ Z(t)^* P Z(t) + \int_{t-D}^{t} Z(s)^* P Z(s) \diff s \right\} \\
& - \dfrac{\gamma_2 (1-\beta) - C_5}{\lambda_M(P)} Z(t-D)^* P Z(t-D) \\
& - \dfrac{\alpha}{2} \sum\limits_{k \geq N_0+1} \left\vert \left< X(t) - B u_D(t) , \psi_k \right>_\mathcal{H} \right\vert^2 \\
& + \dfrac{1}{m_R} \left( \dfrac{2 (m_R + \Vert BK \Vert^2)}{\alpha m_R} + \dfrac{( \gamma_1 (1+D) + \gamma_2 ) \Vert P \Vert^2}{\beta} \right) \\
& \phantom{+}\, \times \underset{\tau \in [t-D,t]}{\mathrm{sup}} \Vert d(\tau) \Vert_\mathcal{H}^2 \\
\leq & - 2 \kappa_0 V(t) + C_6 \underset{\tau \in [t-D,t]}{\mathrm{sup}} \Vert d(\tau) \Vert_\mathcal{H}^2 ,
\end{align*}
where $\lambda_M(P) > 0$ stands for the largest eigenvalue of $P$,
\begin{equation*}
\kappa_0 \triangleq \dfrac{1}{2} \min \left(  
\dfrac{1-\beta}{\lambda_M(P)} , \dfrac{1-\beta  - C_5/\gamma_2}{\lambda_M(P)} , \dfrac{\alpha}{2}
\right)
> 0 ,
\end{equation*}
and
\begin{equation*}
C_6 \triangleq \dfrac{1}{m_R} \left( \dfrac{2 (m_R + \Vert BK \Vert^2)}{\alpha m_R} + \dfrac{( \gamma_1 (1+D) + \gamma_2 ) \Vert P \Vert^2}{\beta} \right) .
\end{equation*}
Then, for all $t > D+t_0$,
\begin{equation}\label{eq: V exp convergence 0}
\dfrac{\mathrm{d}}{\mathrm{d}t}\left[ e^{2 \kappa_0 (\cdot)} V \right](t) \leq C_6 e^{2 \kappa_0 t} \underset{\tau \in [t-D,t]}{\mathrm{sup}} \Vert d(\tau) \Vert_\mathcal{H}^2 .
\end{equation}
As $V \in \mathcal{C}^1(\mathbb{R}_+;\mathbb{R})$, this yields, for all $t \geq D+t_0$,
\begin{align*}
& e^{2 \kappa_0 t} V(t) - e^{2 \kappa_0 (D+t_0)} V(D+t_0) \\
& \qquad\leq C_6 \int_{D+t_0}^{t} e^{2 \kappa_0 s} \underset{\tau \in [s-D,s]}{\mathrm{sup}} \Vert d(\tau) \Vert_\mathcal{H}^2 \diff s \\
& \qquad\leq \dfrac{C_6}{2 \kappa_0} e^{2 \kappa_0 t} \underset{\tau \in [0,t]}{\mathrm{sup}} \Vert d(\tau) \Vert_\mathcal{H}^2 .
\end{align*}
We deduce that, for all $t \geq D+t_0$, 
\begin{equation}\label{eq: V exp convergence 1}
V(t) \leq e^{-2 \kappa_0 (t-D-t_0)} V(D+t_0) + \dfrac{C_6}{2 \kappa_0} \underset{\tau \in [0,t]}{\mathrm{sup}} \Vert d(\tau) \Vert_\mathcal{H}^2 ,
\end{equation}
and thus, from (\ref{eq: X bounded above by cste times sqrt V}) and using the inequality $\sqrt{a+b} \leq \sqrt{a} + \sqrt{b}$ for all $a,b \geq 0$, we obtain that the claimed estimate (\ref{eq: trajectories converge exponentially to zero}) holds true for all $t \geq D+t_0$. Finally, from (\ref{eq: lower bound V}), the control input is such that, for all $t \geq 0$,
\begin{equation}\label{eq: estimation of u from V}
\Vert u(t) \Vert
\leq \Vert K \Vert \Vert Z(t) \Vert 
\leq \dfrac{\Vert K \Vert}{\sqrt{C_2(\gamma_1)}} \sqrt{V(t)} ,
\end{equation}
from which we can deduce that the estimate (\ref{eq: command effort converge exponentially to zero}) is also satisfied for all $t \geq D+t_0$. \qed

\begin{remark}
Coefficient $\beta \in (0,1)$ represents a trade-off between the guaranteed decay rate $\kappa_0$ and the coefficient $C_6 / (2 \kappa_0)$ that reflects the impact of the external disturbance on the system trajectory. In particular, taking $\beta \rightarrow 0^+$ will result in an increasing of the decay rate $\kappa_0$ but also $C_6 / (2 \kappa_0) \rightarrow + \infty$.
\end{remark}

\subsubsection{ISS estimate}
In order to complete the proof of Theorem~\ref{thm: ISS}, we resort to the following lemma that provides an estimate of $V(t)$ over the time interval $[0,d+t_0]$.

\begin{lemma}
Under the assumptions of Theorem~\ref{thm: ISS}, there exist constants $C_9=C_9(\gamma_1,\gamma_2) > 0$ and $C_{10}=C_{10}(\gamma_1,\gamma_2) > 0$, independent of $X_0$ and $d$, such that for all $t \in [0,D+t_0]$, 
\begin{equation}\label{eq: estimate V(t) over [0,D+t0]}
V(t) \leq C_9 \Vert X_0 \Vert^2 + C_{10} \underset{\tau \in [0,t]}{\mathrm{sup}} \Vert d(\tau) \Vert_\mathcal{H}^2 .
\end{equation}
\end{lemma}

\textbf{Proof.} 
With $W(t) = \dfrac{1}{2} \Vert Z(t) \Vert^2$, we have for all $t \geq 0$,
\begin{align*}
\dot{W}(t) 
& = \operatorname{Re} \left< \dot{Z}(t) , Z(t) \right> \\
& \leq \Vert \dot{Z}(t) \Vert \times \Vert Z(t) \Vert \\
& \overset{(\ref{eq: EDO satisfied by Z})}{\leq} \Vert A_{N_0} + \varphi(t) e^{-D A_{N_0}} B_{N_0} K \Vert \times \Vert Z(t) \Vert^2 \\
& \phantom{\leq}\, + \Vert D_{N_0}(t) \Vert \times \Vert Z(t) \Vert \\
& \overset{\mathrm{Y.I.}}{\leq} \left\{ \Vert A_{N_0} \Vert + 1 \times \Vert e^{-D A_{N_0}} B_{N_0} K \Vert \right\} \Vert Z(t) \Vert^2 \\
& \phantom{\overset{\mathrm{Y.I.}}{\leq}}\; + \dfrac{1}{2} \Vert D_{N_0}(t) \Vert^2 + \dfrac{1}{2} \Vert Z(t) \Vert^2 \\
& \overset{(\ref{eq: above estimate of D_N0})}{\leq} 2 C_7 W(t) + \dfrac{1}{2 m_R} \Vert d(t) \Vert_\mathcal{H}^2 
\end{align*}
with $C_7 \triangleq \Vert A_{N_0} \Vert + \Vert e^{-D A_{N_0}} B_{N_0} K \Vert + 1/2 > 0$. Then, for all $t \geq 0$, 
\begin{equation*}
W(t) \leq e^{2 C_7 t} W(0) + \dfrac{1}{4 m_R C_7} e^{2 C_7 t} \underset{\tau \in [0,t]}{\mathrm{sup}} \Vert d(\tau) \Vert_\mathcal{H}^2  .
\end{equation*}
Using (\ref{eq: Riesz basis - inequality}), and from (\ref{eq: Artstein transform}) $Z(0) = Y(0)$, we have $\Vert Z(0) \Vert = \Vert Y(0) \Vert \leq \Vert X_0 \Vert_\mathcal{H} / \sqrt{m_R}$. We deduce that, for all $t \geq 0$, 
\begin{equation}\label{eq: upper bound Z(t)}
\Vert Z(t) \Vert^2 \leq \dfrac{e^{2 C_7 t}}{m_R} \Vert X_0 \Vert^2 + \dfrac{1}{2 m_R C_7} e^{2 C_7 t} \underset{\tau \in [0,t]}{\mathrm{sup}} \Vert d(\tau) \Vert_\mathcal{H}^2 .
\end{equation}

From $u_D(t) = u(t-D) = \varphi(t-D) K Z(t-D)$, we infer that, for all $t \in [0,D+t_0]$,
\begin{equation}\label{eq: upper bound u_D(t)}
\Vert u_D(t) \Vert 
\leq \dfrac{\Vert K \Vert e^{C_7 t_0}}{\sqrt{m_R}} \Vert X_0 \Vert_\mathcal{H} 
+ \dfrac{\Vert K \Vert e^{C_7 t_0}}{\sqrt{2 m_R C_7}} \underset{\tau \in [0,t]}{\mathrm{sup}} \Vert d(\tau) \Vert_\mathcal{H} 
\end{equation}
and, from
\begin{align*}
\dot{u}_D(t) 
& = \dot{\varphi}(t-D) K Z(t-D) + \varphi(t-D) K \dot{Z}(t-D) \\
& = \varphi(t-D) K \left( A_{N_0} + \varphi(t-D) e^{-D A_{N_0}} B_{N_0} K \right) Z(t-D) \\
& \phantom{=}\, + \dot{\varphi}(t-D) K Z(t-D) + \varphi(t-D) K D_{N_0}(t-D) ,
\end{align*}
we obtain that, for all $t \in [0,D+t_0]$,
\begin{align}
\Vert \dot{u}_D(t) \Vert 
& \leq \dfrac{C_8 e^{C_7 t_0}}{\sqrt{m_R}} \Vert X_0 \Vert_\mathcal{H} \nonumber \\
& \phantom{\leq}\, + \dfrac{1}{\sqrt{m_R}} \left( \Vert K \Vert + \dfrac{C_8}{\sqrt{2 C_7}} e^{C_7 t_0} \right) \underset{\tau \in [0,t]}{\mathrm{sup}} \Vert d(\tau) \Vert_\mathcal{H} \label{eq: upper bound dot_u_D(t)}
\end{align}
with $C_8 \triangleq \Vert \dot{\varphi} \Vert_\infty \Vert K \Vert + \Vert K \Vert \left( \Vert A_{N_0} \Vert + \Vert e^{-D A_{N_0}} B_{N_0} K \Vert \right)$.

To conclude, it is sufficient to note that from (\ref{Lyapunov function}), we have for all $t \geq 0$,
\begin{align*}
V(t) 
& \leq \gamma_1 \lambda_M(P) \left\{ \Vert Z(t) \Vert^2 + \int_{t-D}^{t} \varphi(s) \Vert Z(s) \Vert^2 \diff s  \right\} \\
& \phantom{\leq}\, + \gamma_2 \lambda_M(P) \varphi(t-D) \Vert Z(t-D) \Vert^2 \\
& \phantom{\leq}\,  + \dfrac{1}{m_R} \Vert X(t) \Vert_\mathcal{H}^2 
+ \dfrac{\Vert B \Vert^2}{m_R} \Vert u_D(t) \Vert_\mathcal{H}^2 ,
\end{align*}
where, as $X$ is a classical solution of (\ref{def: boundary control system - closed-loop}) and noting that $u_D(0)=u(-D)=0$, we have 
\begin{align*}
X(t) 
& = S(t) X_0 + B u_D(t) \\
& \phantom{=}\, + \int_0^t S(t-\tau) \left\{ -B \dot{u}_D(\tau) + \mathcal{A}B u_D(\tau) + d(\tau) \right\} \diff \tau .
\end{align*}
By direct estimation and using (\ref{eq: upper bound Z(t)}-\ref{eq: upper bound dot_u_D(t)}), we deduce that the conclusion of the lemma holds true. \qed

We can now complete the proof of Theorem~\ref{thm: ISS}. Indeed, for a given arbitrary $\beta \in (0,1)$ and by selecting $\gamma_1 > C_1/\lambda_m(P)$ and $\gamma_2 > \max \left( \Vert B K \Vert^2/(m_R \lambda_m(P)) , C_5/(1-\beta) \right)$, we obtain from (\ref{eq: V exp convergence 1}) and (\ref{eq: estimate V(t) over [0,D+t0]}) that the following estimate holds true
\begin{equation*}
V(t) \leq C_9 e^{-2 \kappa_0 (t - D - t_0)} \Vert X_0 \Vert_\mathcal{H}^2 
+ \left( \dfrac{C_6}{2 \kappa_0} + C_{10} \right) \underset{\tau \in [0,t]}{\mathrm{sup}} \Vert d(\tau) \Vert_\mathcal{H}^2,
\end{equation*}
for all $t \geq 0$. From (\ref{eq: X bounded above by cste times sqrt V}), we obtain that, for all $t \geq 0$, 
\begin{align*}
\left\Vert X(t) \right\Vert_\mathcal{H} 
& \leq 
\left\{ C_4 \sqrt{C_9} e^{\kappa_0 (D + t_0)} \right\} e^{- \kappa_0 t} \Vert X_0 \Vert_\mathcal{H} \\
& \phantom{\leq}\, + C_4 \sqrt{\dfrac{C_6}{2 \kappa_0} + C_{10}} \underset{\tau \in [0,t]}{\mathrm{sup}} \Vert d(\tau) \Vert_\mathcal{H} .
\end{align*}
It shows that the claimed ISS estimate (\ref{eq: global exponential stability}) holds true. The estimate of the control input (\ref{eq: command effort global exponential convergence}) follows from (\ref{eq: estimation of u from V}), which concludes the proof of Theorem~\ref{thm: ISS}.

\section{Application to the stability analysis of a closed-loop interconected IDS-ODE system}\label{sec: Application to the stability analysis of a closed-loop interconected IDS-ODE system}
As an application of the ISS property of the closed-loop system (\ref{def: boundary control system - closed-loop}), we propose to study the stability of a related IDS-ODE interconnection. Specifically, we consider the case where the external input $d$ depends on the state of an ODE satisfying a certain ISS estimate.

\subsection{Dynamics of the closed-loop interconnected IDS-ODE system and well-posedness}
Let $D,t_0 > 0$ be given. We consider a given transition signal $\varphi \in \mathcal{C}^2([-D,+\infty);\mathbb{R})$ such that $0 \leq \varphi \leq 1$, $\left. \varphi \right\vert_{[-D,0]} = 0$, and $\left. \varphi \right\vert_{[t_0,+\infty)} = 1$. Let $f_1 \in \mathcal{C}^1(\mathbb{K}^n \times \mathcal{H} \times \mathbb{K}^{m_v} ; \mathbb{K}^n)$ and $f_2 \in \mathcal{C}^1(\mathbb{K}^n \times \mathcal{H} \times \mathbb{K}^{m_v} , \mathcal{H})$ be two vector fields. We make the following assumption.

\begin{assum}\label{assum: A4}
\begin{enumerate}
\item Vector fields $f_1(x,X,v)$ and $f_2(x,X,v)$ are (globally) Lipschitz continuous in $(x,X)$ on $\mathbb{K}^n \times \mathcal{H}$, uniformly in $v$ over any compact subset of $\mathbb{K}^{m_v}$.
\item There exist constants $D_1,D_2,D_3 \geq 0$ such that, for all $x \in \mathbb{K}^n$, $ X \in \mathcal{H}$, and $v \in \mathbb{K}^{m_v}$,
\begin{equation}\label{eq: HYP - f2}
\Vert f_2(x,X,v) \Vert_\mathcal{H} 
\leq D_1 \Vert x \Vert + D_2 \Vert X \Vert_\mathcal{H} + D_3 \Vert v \Vert .
\end{equation}
\item The ODE $\dot{x} = f_1(x,X,v)$ is such that there exist $\tilde{\kappa}_0 , \tilde{C}_0 , \tilde{C}_1 , \tilde{C}_2 \in \mathbb{R}_+$ such that, for any given initial condition $x_0 \in \mathbb{K}^n$ and functions $X \in \mathcal{C}^0(\mathbb{R}_+;\mathcal{H})$ and $v \in \mathcal{C}^0(\mathbb{R}_+;\mathbb{K}^{m_v})$, the following ISS estimate holds true for all $t \geq 0$
\end{enumerate}
\begin{equation}\label{eq: HYP - ISS estimate f1}
\Vert x(t) \Vert^2 
\leq \tilde{C}_0^2 e^{-2 \tilde{\kappa}_0 t} \Vert x_0 \Vert^2 
+ \sup\limits_{\tau \in [0,t]} \left\{ \tilde{C}_1^2 \Vert X(\tau) \Vert_\mathcal{H}^2 + \tilde{C}_2^2 \Vert v(\tau) \Vert^2 \right\} .
\end{equation}
\end{assum}

Note that the above assumption implies that $\tilde{C}_0 \geq 1$. The considered closed-loop system takes the following form:
\begin{equation}\label{def: interconnected IDS-ODE - closed-loop}
\left\{\begin{split}
\dot{x}(t) & = f_1(x(t),X(t),v(t)) , \\
\dfrac{\mathrm{d} X}{\mathrm{d} t}(t) & = \mathcal{A} X(t) + f_2(x(t),X(t),v(t)) , \\
\mathcal{B} X (t) & = u_D(t) = u(t-D) , \\
\left. u \right\vert_{[-D,0]} & = 0 \\
u(t) & = \varphi(t) K Y(t) \\
& \phantom{=} + \varphi(t) K \int_{\max(t-D,0)}^{t} e^{(t-s-D)A_{N_0}} B_{N_0} u(s) \diff s , \\
Y(t) & = 
\begin{bmatrix}
\left< X(t) , \psi_1 \right>_\mathcal{H} \\ \vdots \\ \left< X(t) , \psi_{N_0} \right>_\mathcal{H}
\end{bmatrix} , \\
x(0) & = x_0 , \\
X(0) & = X_0 
\end{split}\right.
\end{equation}
for $t \geq 0$. The feedback gain $K \in \mathbb{K}^{m \times N_0}$ is such that $A_\mathrm{cl} \triangleq A_{N_0} + e^{-D A_{N_0}} B_{N_0} K$ is Hurwitz (with desired pole placement). Function $u$ still represents the control input while function $v : \mathbb{R}_+ \rightarrow \mathbb{K}^{m_v}$ represents a disturbance.

The well-posedness of the closed-loop system (\ref{def: interconnected IDS-ODE - closed-loop}) is assessed via the following result.

\begin{lemma}\label{lemma: well-posedness closed-loop IDS-ODE}
Let $(\mathcal{A},\mathcal{B})$ be an abstract boundary control system and $f_1 \in \mathcal{C}^1(\mathbb{K}^n \times \mathcal{H} \times \mathbb{K}^{m_v} ; \mathbb{K}^n)$ and $f_2 \in \mathcal{C}^1(\mathbb{K}^n \times \mathcal{H} \times \mathbb{K}^{m_v} , \mathcal{H})$ be vector fields such that Assumptions~\ref{assum: A1}, \ref{assum: A2}, \ref{assum: A3}, and~\ref{assum: A4} hold true. For any $(x_0,X_0) \in \mathbb{K}^n \times D(\mathcal{A}_0)$ and $v \in \mathcal{C}^1(\mathbb{R}_+;\mathbb{K}^{m_v})$, the closed-loop system (\ref{def: interconnected IDS-ODE - closed-loop}) has a unique classical solution $(x,X) \in \mathcal{C}^1(\mathbb{R}_+;\mathbb{K}^n) \times \left( \mathcal{C}^0(\mathbb{R}_+;D(\mathcal{A})) \cap \mathcal{C}^1(\mathbb{R}_+;\mathcal{H}) \right)$. Introducing  $d(t) \triangleq f_2(x(t),X(t),v(t))$, we have $d \in \mathcal{C}^1(\mathbb{R}_+;\mathcal{H})$. Thus, $X$ is the classical solution of (\ref{def: boundary control system - closed-loop}) associated with the initial condition $X_0$ and the distributed disturbance $d$. Consequently, both Lemma~\ref{lemma: well-posedness closed-loop PDE} and Theorem~\ref{thm: ISS} apply to $X$.
\end{lemma}

The proof of Lemma~\ref{lemma: well-posedness closed-loop IDS-ODE} follows from the same arguments as the one used in the proof of Lemma~\ref{lemma: well-posedness closed-loop PDE} and from classical theorems on the existence and uniqueness of classical solutions for lipschitz perturbations of linear evolution equations, see, e.g., \cite[Th.~1.2 and Th.~1.5, Chap.~6]{pazy2012semigroups}.

\begin{remark}\label{rem: blow up in finite time}
The first point of Assumption~\ref{assum: A4} regarding the Lipschitz continuity of vector fields $f_1,f_2$ is used to ensure the existence of solutions defined over $\mathbb{R}_+$. In particular, it avoids any potential blow up of the solution in finite time. If this assumption is removed, the existence of the classical solution is \emph{a priori} only guaranteed over a time interval $[0,t_{\max})$ with $0 < t_{\max} \leq + \infty$. Furthermore, if $t_{\max} < +\infty$, we have the blow up of the solution in finite time, i.e., $\Vert x(t) \Vert + \Vert X(t) \Vert_\mathcal{H} \underset{t \rightarrow (t_{\max})^{-}}{\longrightarrow} +\infty$, see, e.g., \cite[Th.~1.4 and Th.~1.5, Chap.~6]{pazy2012semigroups}. In this case, the reasoning presented next still applies over the time interval $[0,t_{\max})$ at the condition that no blow up occurs over the time interval $[0,D+t_0]$, i.e., $t_{\max} > D+t_0$. This can be ensured by assuming that the following small gain condition holds true:
\begin{equation}\label{eq: sufficient small gain condition for no blow up before D+t0}
( D_1 \tilde{C}_1 + D_2 ) C_4 \sqrt{C_{10}} < 1  .
\end{equation}
Indeed, from (\ref{eq: X bounded above by cste times sqrt V}), (\ref{eq: estimate V(t) over [0,D+t0]}), and (\ref{eq: HYP - f2}-\ref{eq: HYP - ISS estimate f1}), we obtain that, for all $t \in [0,D+t_0] \cap [0,t_{\max})$,
\begin{align*}
\Vert X(t) \Vert_\mathcal{H} 
& \leq 
D_1 \tilde{C}_0 C_4 \sqrt{C_{10}} \Vert x_0 \Vert 
+ C_4 \sqrt{C_9} \Vert X_0 \Vert_\mathcal{H} \\
& \phantom{\leq}\, + (D_1 \tilde{C}_1 + D_2) C_4 \sqrt{C_{10}} \sup\limits_{\tau \in [0,t]} \Vert X(\tau) \Vert_\mathcal{H}\\
& \phantom{\leq}\, + (D_1 \tilde{C}_2 + D_3) C_4 \sqrt{C_{10}} \sup\limits_{\tau \in [0,t]} \Vert v(\tau) \Vert .
\end{align*}
Under the small gain assumption (\ref{eq: sufficient small gain condition for no blow up before D+t0}), we can introduce
\begin{equation}
\Gamma \triangleq \left( 1 - ( D_1 \tilde{C}_1 + D_2 ) C_4 \sqrt{C_{10}} \right)^{-1} > 0 ,
\end{equation}
which yields
\begin{align*}
& \sup\limits_{\tau \in [0,D+t_0] \cap [0,t_{\max})} \Vert X(\tau) \Vert_\mathcal{H} \\
& \qquad\qquad\leq 
\Gamma D_1 \tilde{C}_0 C_4 \sqrt{C_{10}} \Vert x_0 \Vert 
+ \Gamma C_4 \sqrt{C_9} \Vert X_0 \Vert_\mathcal{H} \\
& \qquad\qquad\phantom{\leq}\, + \Gamma  (D_1 \tilde{C}_2 + D_3) C_4 \sqrt{C_{10}} \sup\limits_{\tau \in [0,D+t_0]} \Vert v(\tau) \Vert \\
& \qquad\qquad< \infty .
\end{align*}
From (\ref{eq: HYP - ISS estimate f1}) we infer that
\begin{equation*}
\sup\limits_{\tau \in [0,D+t_0] \cap [0,t_{\max})} \left\{ \Vert x(\tau) \Vert + \Vert X(\tau) \Vert_\mathcal{H} \right\} < \infty ,
\end{equation*}
and, consequently, $t_{\max} > D+t_0$.
\end{remark}

\subsection{Small gain condition ensuring the stability of the IDS-ODE interconnection}

The objective of this section is to demonstrate the following result.

\begin{theorem}\label{thm: stab IDS-ODE loop}
Let $(\mathcal{A},\mathcal{B})$ be an abstract boundary control system and $f_1 \in \mathcal{C}^1(\mathbb{K}^n \times \mathcal{H} \times \mathbb{K}^{m_v} ; \mathbb{K}^n)$ and $f_2 \in \mathcal{C}^1(\mathbb{K}^n \times \mathcal{H} \times \mathbb{K}^{m_v} , \mathcal{H})$ be vector fields such that Assumptions~\ref{assum: A1}, \ref{assum: A2}, \ref{assum: A3}, and~\ref{assum: A4} hold true. We assume that the small gain condition
\begin{equation}\label{eq: sufficient small gain condition}
( D_1 \tilde{C}_1 + D_2 ) C_4 \sqrt{\dfrac{C_6}{2 \kappa_0}} < 1 
\end{equation}
is satisfied. Then, there exist constants $\delta_\epsilon \in (0,\kappa_0)$ and $G_i,H_i \in \mathbb{R}_+$, $0 \leq i \leq 3$, such that, for any $(x_0,X_0) \in \mathbb{K}^n \times D(\mathcal{A}_0)$ and $v \in \mathcal{C}^1(\mathbb{R}_+;\mathbb{K}^{m_v})$, the classical solution $(x,X)$ of (\ref{def: interconnected IDS-ODE - closed-loop}) associated with the initial condition $(x_0,X_0)$ and the disturbance $v$ satisfies for all $t \geq D+t_0$ the following fading memory estimate:
\begin{align}
\Vert x(t) \Vert + \Vert X(t) \Vert_\mathcal{H}
& \leq G_0 e^{-\delta_\epsilon t} ( \Vert x_0 \Vert + \Vert X_0 \Vert_\mathcal{H}) \nonumber \\
& \phantom{\leq}\, + G_1 e^{-\delta_\epsilon t} \sup\limits_{\tau \in [0,D+t_0]} \Vert x(\tau) \Vert \label{eq: IDS-ODE stab estimate 1}  \\
& \phantom{\leq}\, + G_2 e^{-\delta_\epsilon t} \sup\limits_{\tau \in [0,D+t_0]} \Vert X(\tau) \Vert_\mathcal{H} \nonumber  \\
& \phantom{\leq}\, + G_3 \sup\limits_{\tau \in [0,t]} e^{- \delta_\epsilon ( t - \tau )} \Vert v(\tau) \Vert , \nonumber
\end{align}
and the control law satisfies
\begin{align}
\Vert u(t) \Vert
& \leq H_0 e^{-\delta_\epsilon t} ( \Vert x_0 \Vert + \Vert X_0 \Vert_\mathcal{H}) \nonumber \\
& \phantom{\leq}\, + H_1 e^{-\delta_\epsilon t} \sup\limits_{\tau \in [0,D+t_0]} \Vert x(\tau) \Vert \label{eq: IDS-ODE stab estimate 2} \\
& \phantom{\leq}\, + H_2 e^{-\delta_\epsilon t} \sup\limits_{\tau \in [0,D+t_0]} \Vert X(\tau) \Vert_\mathcal{H} \nonumber \\
& \phantom{\leq}\, + H_3 \sup\limits_{\tau \in [0,t]} e^{- \delta_\epsilon ( t - \tau )} \Vert v(\tau) \Vert \nonumber
\end{align}
for all $t \geq D+t_0$.
\end{theorem}

\begin{remark}
As the system is in open loop over the time interval $[0,D]$ and then the time interval $[D,D+t_0]$ is employed to switch from open loop to closed loop, we can interpret $\left. x \right\vert_{[0,D+t_0]}$ and $\left. X \right\vert_{[0,D+t_0]}$ as initial perturbations. In this case, (\ref{eq: IDS-ODE stab estimate 1}) can be seen as an ISS estimate with fading memory with respect to the initial perturbations $\left. x \right\vert_{[0,D+t_0]}$ and $\left. X \right\vert_{[0,D+t_0]}$ and the disturbance $v$.
\end{remark}

The remaining of this section is devoted to the proof of Theorem~\ref{thm: stab IDS-ODE loop} through an adaptation of the approach presented in~\cite{karafyllis2018input} for the study of the stability of IDS-ODE or PDE-PDE interconnections via a small gain approach. In order to be able to apply the results of the previous section, $V$ is still defined by (\ref{Lyapunov function}) with $\gamma_1 , \gamma_2$ large enough\footnote{More precisely, they are selected such that $\gamma_1 > C_1/\lambda_m(P)$ and $\gamma_2 > \max \left( \Vert B K \Vert^2/(m_R \lambda_m(P)) , C_5/(1-\beta) \right)$.}.

\subsubsection{Conversion of the ISS estimates into fading memory estimates}
Following the methodology presented in~\cite{karafyllis2018input} for studying the stability of IDS-ODE or PDE-PDE interconnections, the key step relies in the conversion of the ISS estimates satisfied by each component of the interconnections into fading memory estimates via the following lemma~\cite[Lemma~7.1]{karafyllis2018input}.

\begin{lemma}[Conversion Lemma]
For every $\sigma > 0$, $M \geq 1$, and $\epsilon > 0$, there exists a constant $\delta \in (0,\sigma)$ such that for any continuous functions $\phi : \mathbb{R}_+ \rightarrow \mathbb{R}_+$ and $y : \mathbb{R}_+ \rightarrow \mathbb{R}_+$ for which there exists a constant $\gamma \geq 0$ such that the following inequality holds true for all $t_0 \geq 0$ and $t \geq t_0$,
\begin{equation}\label{eq: conversion lemma condition application}
\phi(t) \leq M e^{- \sigma (t-t_0)} \phi(t_0) + \gamma \sup\limits_{t_0 \leq s \leq t} y(s) ,
\end{equation}
then the following inequality holds for all $t \geq t_0$:
\begin{equation*}
\phi(t) \leq M e^{- \delta t} \phi(0) + \gamma (1+\epsilon) \sup\limits_{0 \leq s \leq t} e^{- \delta (t-s)} y(s) .
\end{equation*}
\end{lemma}

Even if the trajectories $X$ of (\ref{def: boundary control system - closed-loop}) satisfy the ISS estimate (\ref{eq: global exponential stability}) provided by Theorem~\ref{thm: ISS}, we cannot directly apply the Conversion Lemma because the semigroup property does not hold true. This is due to the time-varying nature of (\ref{def: boundary control system - closed-loop}) induced by the transition from open loop to closed loop via $\varphi$, yielding $\left. u_D \right\vert_{[0,D)} = 0$. Therefore, we cannot directly deduce from the ISS estimate (\ref{eq: global exponential stability}) that an estimate similar to (\ref{eq: conversion lemma condition application}) holds true for all $t \geq t_0 \geq 0$. In order to avoid this pitfall, we are not going to apply the Conversion Lemma to the system trajectories $X$ but to the Lyapunov function $V$. Indeed, with $d(t) = f_2(x(t),X(t),v(t))$, we know from Lemma~\ref{lemma: well-posedness closed-loop IDS-ODE} that $X$ is solution of (\ref{def: boundary control system - closed-loop}) associated with the initial condition $X_0$ and the distributed disturbance $d$. Consequently, we deduce from (\ref{eq: V exp convergence 0}) that, for all $t_2 \geq t_1 \geq D+t_0$,
\begin{align*}
& e^{2 \kappa_0 t_2} V(t_2) - e^{2 \kappa_0 t_1} V(t_1) \\
& \leq C_6 \int_{t_1}^{t_2} e^{2 \kappa_0 s} \sup\limits_{\tau \in [s-D,s]} \Vert d(\tau) \Vert_\mathcal{H}^2 \diff s \\
& \leq \dfrac{C_6}{2 \kappa_0} e^{2 \kappa_0 t_2} \sup\limits_{s \in [t_1,t_2]} \, \sup\limits_{\tau \in [s-D,s]} \Vert d(\tau) \Vert_\mathcal{H}^2 .
\end{align*}
This yields, for all $t_2 \geq t_1 \geq D+t_0$,
\begin{equation*}
V(t_2) 
\leq e^{- 2 \kappa_0 (t_2 - t_1)} V(t_1) 
+ \dfrac{C_6}{2 \kappa_0} \sup\limits_{s \in [t_1,t_2]} \, \sup\limits_{\tau \in [s-D,s]} \Vert d(\tau) \Vert_\mathcal{H}^2 .
\end{equation*}
Introducing $\hat{\kappa}_0 = \min ( \kappa_0 , \tilde{\kappa}_0 ) > 0$ and noting that $\tilde{C}_0 \geq 1$, then we have for all $t_2 \geq t_1 \geq 0$,
\begin{align}
V(t_2 + (D+t_0)) 
& \leq \tilde{C}_0^2 e^{- 2 \hat{\kappa}_0 (t_2 - t_1)} V(t_1 + (D+t_0)) \label{eq: inequality 1 to apply conv lemma} \\
& \phantom{\leq} + \dfrac{C_6}{2 \kappa_0} \sup\limits_{s \in [t_1,t_2]} \, \sup\limits_{\tau \in [s+t_0,s+(D+t_0)]} \Vert d(\tau) \Vert_\mathcal{H}^2 . \nonumber
\end{align}
Furthermore, as the trajectories of the ODE $\dot{x} = f_1(x,X,v)$ satisfy the semigroup property, we also have from (\ref{eq: HYP - ISS estimate f1}) that\footnote{We estimate by replacing $\tilde{\kappa}_0$ by $\hat{\kappa}_0$.} for all $t_2 \geq t_1 \geq 0$,
\begin{align}
\Vert x(t_2) \Vert^2 
& \leq \tilde{C}_0^2 e^{-2 \hat{\kappa}_0 (t_2 - t_1)} \Vert x(t_1) \Vert^2 \label{eq: inequality 2 to apply conv lemma} \\
& \phantom{\leq}\, + \sup\limits_{\tau \in [t_1,t_2]} \left\{ \tilde{C}_1^2 \Vert X(\tau) \Vert_\mathcal{H}^2 + \tilde{C}_2^2 \Vert v(\tau) \Vert^2 \right\} . \nonumber
\end{align}

\begin{remark}
The introduction of the constant $\tilde{C}_0^2 \geq 1$ in (\ref{eq: inequality 1 to apply conv lemma}) is motivated by the will to apply the Conversion Lemma simultaneously to both (\ref{eq: inequality 1 to apply conv lemma}-\ref{eq: inequality 2 to apply conv lemma}). Even if this yields some conservatism is the estimate with respect to the value of $V$ at the lower bound of the interval of integration, such an introduction will have no impact on the conservatism of the small gain condition (\ref{eq: sufficient small gain condition}).
\end{remark}

We now apply the Conversion Lemma. For $\sigma = 2 \hat{\kappa}_0$ and $M = \tilde{C}_0^2 \geq 1$, we denote by $2 \delta_\epsilon \in (0 , 2 \hat{\kappa}_0)$ the constant ``$\delta$'' provided by the Conversion Lemma (which is independent of $x_0$, $X_0$, and $v$) for any given $\epsilon > 0$.  
From the proof of the Conversion Lemma in~\cite[Lemma~7.1]{karafyllis2018input}, we can select $\delta_\epsilon$ such that $\delta_\epsilon \underset{\epsilon \rightarrow 0^+}{\longrightarrow} 0^+$. 

Applying the Conversion Lemma to (\ref{eq: inequality 1 to apply conv lemma}) with $\phi(t) = V(t + (D+t_0))$, $y(t) = \sup\limits_{\tau \in [t+t_0,t+(D+t_0)]} \Vert d(\tau) \Vert_\mathcal{H}^2$, and $\gamma = C_6 / (2 \kappa_0)$, we infer that, for all $t \geq 0$, 
\begin{align*}
& e^{2 \delta_\epsilon t} V(t + (D+t_0)) \\
& \leq \tilde{C}_0^2 V(D+t_0) \\
& \phantom{\leq}\, + \dfrac{C_6}{2 \kappa_0}(1+\epsilon) \sup\limits_{s \in [0,t]} \left\{ e^{2 \delta_\epsilon s} \sup\limits_{\tau \in [s+t_0,s+(D+t_0)]} \Vert d(\tau) \Vert_\mathcal{H}^2 \right\}  .
\end{align*}
Noting that $s + t_0 \leq \tau$ implies $s \leq \tau - t_0$ and thus $e^{2 \delta_\epsilon s} \leq e^{2 \delta_\epsilon \tau} e^{- 2 \delta_\epsilon t_0}$, we obtain for all $t \geq 0$,
\begin{align}
& e^{2 \delta_\epsilon t} V(t + (D+t_0)) \nonumber \\
& \leq \tilde{C}_0^2 V(D+t_0) 
+ \dfrac{C_6}{2 \kappa_0}(1+\epsilon) e^{- 2 \delta_\epsilon t_0} \sup\limits_{\tau \in [t_0,t+(D+t_0)]} e^{2 \delta_\epsilon \tau} \Vert d(\tau) \Vert_\mathcal{H}^2  . \label{eq: conv lemma - applied to V}
\end{align}
From (\ref{eq: X bounded above by cste times sqrt V}) we obtain that, for all $t \geq 0$,
\begin{align}
& e^{\delta_\epsilon t} \left\Vert X(t + (D+t_0)) \right\Vert_\mathcal{H}  \nonumber \\
& \leq C_4 \tilde{C}_0 \sqrt{V(D+t_0)} \label{eq: conv lemma - res1} \\
& \phantom{\leq}\, + C_4 \sqrt{\dfrac{C_6}{2 \kappa_0} (1+\epsilon)} e^{-\delta_\epsilon t_0} \sup\limits_{\tau \in [t_0,t+(D+t_0)]} e^{\delta_\epsilon \tau} \Vert d(\tau) \Vert_\mathcal{H} . \nonumber
\end{align}

From the application of the Conversion Lemma to (\ref{eq: inequality 2 to apply conv lemma}) with $\phi(t) = \Vert x(t) \Vert^2$, $y(t) = \tilde{C}_1^2 \Vert X(t) \Vert_\mathcal{H}^2 + \tilde{C}_2^2 \Vert v(t) \Vert^2$, and $\gamma = 1$, we infer that, for all $t \geq 0$,
\begin{align*}
& e^{2 \delta_\epsilon t}  \Vert x(t) \Vert^2 \\ 
& \leq \tilde{C}_0^2 \Vert x_0 \Vert^2
+ (1+\epsilon) \sup\limits_{\tau \in [0,t]} e^{2 \delta_\epsilon \tau} \left\{ \tilde{C}_1^2 \Vert X(\tau) \Vert_\mathcal{H}^2 + \tilde{C}_2^2 \Vert v(\tau) \Vert^2 \right\} .
\end{align*}
This yields, for all $t \geq 0$,
\begin{align}
e^{\delta_\epsilon t} \Vert x(t) \Vert 
& \leq \tilde{C}_0 \Vert x_0 \Vert 
+ \tilde{C}_1 \sqrt{1+\epsilon} \sup\limits_{\tau \in [0,t]} e^{\delta_\epsilon \tau} \Vert X(\tau) \Vert_\mathcal{H} \label{eq: conv lemma - res2} \\
& \phantom{\leq}\, + \tilde{C}_2 \sqrt{1+\epsilon} \sup\limits_{\tau \in [0,t]} e^{\delta_\epsilon \tau} \Vert v(\tau) \Vert . \nonumber
\end{align}

\subsubsection{Stability of the interconnected IDS-ODE}
We can now proceed to the proof of Theorem~\ref{thm: stab IDS-ODE loop}. From (\ref{eq: HYP - f2}) and (\ref{eq: conv lemma - res2}) we obtain that, for all $t \geq 0$,
\begin{align}
& e^{\delta_\epsilon t} \Vert d(t) \Vert_\mathcal{H} \nonumber \\
& = e^{\delta_\epsilon t} \Vert f_2(x(t),X(t),v(t)) \Vert_\mathcal{H} \nonumber \\
& \leq D_1 e^{\delta_\epsilon t} \Vert x(t) \Vert 
+ D_2 e^{\delta_\epsilon t} \Vert X(t) \Vert_\mathcal{H}
+ D_3 e^{\delta_\epsilon t} \Vert v(t) \Vert \nonumber \\
& \leq D_1 \tilde{C}_0 \Vert x_0 \Vert 
+ D_1 \tilde{C}_1 \sqrt{1+\epsilon} \sup\limits_{\tau \in [0,t]} e^{\delta_\epsilon \tau} \Vert X(\tau) \Vert_\mathcal{H} \nonumber \\
& \phantom{\leq}\, + D_1 \tilde{C}_2 \sqrt{1+\epsilon} \sup\limits_{\tau \in [0,t]} e^{\delta_\epsilon \tau} \Vert v(\tau) \Vert
+ D_2 e^{\delta_\epsilon t} \Vert X(t) \Vert_\mathcal{H} \nonumber \\
& \phantom{\leq}\, + D_3 e^{\delta_\epsilon t} \Vert v(t) \Vert \nonumber \\
& \leq D_1 \tilde{C}_0 \Vert x_0 \Vert 
+ ( D_1 \tilde{C}_1 \sqrt{1+\epsilon} + D_2 ) \sup\limits_{\tau \in [0,t]} e^{\delta_\epsilon \tau} \Vert X(\tau) \Vert_\mathcal{H}  \nonumber \\
& \phantom{\leq} + ( D_1 \tilde{C}_2 \sqrt{1+\epsilon} + D_3 ) \sup\limits_{\tau \in [0,t]} e^{\delta_\epsilon \tau} \Vert v(\tau) \Vert . \label{eq: loop IDS-ODE estimation d}
\end{align}
This yields, for all $t \geq 0$,
\begin{align*}
& \sup\limits_{\tau \in [t_0,t+(D+t_0)]} e^{\delta_\epsilon \tau} \Vert d(\tau) \Vert_\mathcal{H} \\
& \leq D_1 \tilde{C}_0 \Vert x_0 \Vert 
+ ( D_1 \tilde{C}_1 \sqrt{1+\epsilon} + D_2) \sup\limits_{\tau \in [0,t+(D+t_0)]} e^{\delta_\epsilon \tau} \Vert X(\tau) \Vert_\mathcal{H} \\
& \phantom{\leq}\, + ( D_1 \tilde{C}_2 \sqrt{1+\epsilon} + D_3 ) \sup\limits_{\tau \in [0,t+(D+t_0)]} e^{\delta_\epsilon \tau} \Vert v(\tau) \Vert \\ 
& \leq D_1 \tilde{C}_0 \Vert x_0 \Vert 
+ ( D_1 \tilde{C}_1 \sqrt{1+\epsilon} + D_2) \sup\limits_{\tau \in [0,D+t_0]} e^{\delta_\epsilon \tau} \Vert X(\tau) \Vert_\mathcal{H} \\
& \phantom{\leq}\, + ( D_1 \tilde{C}_1 \sqrt{1+\epsilon} + D_2 ) \sup\limits_{\tau \in [D+t_0,t+(D+t_0)]} e^{\delta_\epsilon \tau} \Vert X(\tau) \Vert_\mathcal{H} \\
& \phantom{\leq}\, + ( D_1 \tilde{C}_2 \sqrt{1+\epsilon} + D_3 ) \sup\limits_{\tau \in [0,t+(D+t_0)]} e^{\delta_\epsilon \tau} \Vert v(\tau) \Vert .
\end{align*}
Therefore, we deduce from (\ref{eq: conv lemma - res1}) that, for all $t \geq 0$,
\begin{align*}
& \sup\limits_{\tau \in [D+t_0,t+(D+t_0)]} e^{\delta_\epsilon \tau} \Vert X(\tau) \Vert_\mathcal{H} \\
& \leq C_4 \tilde{C}_0 e^{\delta_\epsilon (D+t_0)} \sqrt{V(D+t_0)}
+ D_1 \tilde{C}_0 C_4 \sqrt{\dfrac{C_6}{2 \kappa_0} (1+\epsilon)} e^{\delta_\epsilon D} \Vert x_0 \Vert \\
& \phantom{\leq}\, + ( D_1 \tilde{C}_1 \sqrt{1+\epsilon} + D_2 ) C_4 \sqrt{\dfrac{C_6}{2 \kappa_0} (1+\epsilon)} e^{\delta_\epsilon D} \\ 
& \phantom{\leq \, +}\, \times \sup\limits_{\tau \in [0,D+t_0]} e^{\delta_\epsilon \tau} \Vert X(\tau) \Vert_\mathcal{H} \\
& \phantom{\leq}\, + ( D_1 \tilde{C}_1 \sqrt{1+\epsilon} + D_2 ) C_4 \sqrt{\dfrac{C_6}{2 \kappa_0} (1+\epsilon)} e^{\delta_\epsilon D} \\ 
& \phantom{\leq \, +}\, \times \sup\limits_{\tau \in [D+t_0,t+(D+t_0)]} e^{\delta_\epsilon \tau} \Vert X(\tau) \Vert_\mathcal{H} \\
& \phantom{\leq}\, + ( D_1 \tilde{C}_2 \sqrt{1+\epsilon} + D_3 ) C_4 \sqrt{\dfrac{C_6}{2 \kappa_0} (1+\epsilon)} e^{\delta_\epsilon D} \\ 
& \phantom{\leq \, +}\, \times \sup\limits_{\tau \in [0,t+(D+t_0)]} e^{\delta_\epsilon \tau} \Vert v(\tau) \Vert .
\end{align*}
As $\delta_\epsilon \underset{\epsilon \rightarrow 0^+}{\longrightarrow} 0^+$ and because of the small gain assumption (\ref{eq: sufficient small gain condition}), there exist $\epsilon > 0$ such that 
\begin{equation*}
( D_1 \tilde{C}_1 \sqrt{1+\epsilon} + D_2 ) C_4 \sqrt{\dfrac{C_6}{2 \kappa_0} (1+\epsilon)} e^{\delta_\epsilon D} < 1 .
\end{equation*}
We fix such $\epsilon > 0$, which is independent of the initial condition $(x_0,X_0)$ and the disturbance $v$. Therefore, we obtain that, for all $t \geq 0$, 
\begin{align}
& e^{\delta_\epsilon (t+(D+t_0))} \Vert X(t+(D+t_0)) \Vert_\mathcal{H} \nonumber \\
& \leq \sup\limits_{\tau \in [D+t_0,t+(D+t_0)]} e^{\delta_\epsilon \tau} \Vert X(\tau) \Vert_\mathcal{H} \nonumber \\
& \leq E_1 \sqrt{V(D+t_0)}
+ E_2 \Vert x_0 \Vert 
+ E_3 \sup\limits_{\tau \in [0,D+t_0]} e^{\delta_\epsilon \tau} \Vert X(\tau) \Vert_\mathcal{H} \nonumber \\
& \phantom{\leq}\, + E_4 \sup\limits_{\tau \in [0,t+(D+t_0)]} e^{\delta_\epsilon \tau} \Vert v(\tau) \Vert .  \label{eq: IDS-ODE stab estimate 1 - intermediate 1} 
\end{align}
where 
\begin{align*}
E_1 & = \Delta C_4 \tilde{C}_0 e^{\delta_\epsilon (D+t_0)}, \\
E_2 & = \Delta D_1 \tilde{C}_0 C_4 \sqrt{\dfrac{C_6}{2 \kappa_0} (1+\epsilon)} e^{\delta_\epsilon D} , \\
E_3 & = \Delta ( D_1 \tilde{C}_1 \sqrt{1+\epsilon} + D_2 ) C_4 \sqrt{\dfrac{C_6}{2 \kappa_0} (1+\epsilon)} e^{\delta_\epsilon D} , \\
E_4 & = \Delta ( D_1 \tilde{C}_2 \sqrt{1+\epsilon} + D_2 ) C_4 \sqrt{\dfrac{C_6}{2 \kappa_0} (1+\epsilon)} e^{\delta_\epsilon D} ,
\end{align*}
with $\Delta > 0$ defined by
\begin{equation*}
\Delta = \left( 1 - ( D_1 \tilde{C}_1 \sqrt{1+\epsilon} + D_2 ) C_4 \sqrt{\dfrac{C_6}{2 \kappa_0} (1+\epsilon)} e^{\delta_\epsilon D} \right)^{-1} .
\end{equation*}
From (\ref{eq: conv lemma - res2}), we have, for all $t \geq 0$, 
\begin{align}
& e^{\delta_\epsilon (t+(D+t_0))} \Vert x(t+(D+t_0)) \Vert \nonumber \\
& \leq \tilde{C}_0 \Vert x_0 \Vert 
+ \tilde{C}_1 \sqrt{1+\epsilon} \sup\limits_{\tau \in [0,D+t_0]} e^{\delta_\epsilon \tau} \Vert X(\tau) \Vert_\mathcal{H} \nonumber \\
& \phantom{\leq}\, + \tilde{C}_1 \sqrt{1+\epsilon} \sup\limits_{\tau \in [D+t_0,t+(D+t_0)]} e^{\delta_\epsilon \tau} \Vert X(\tau) \Vert_\mathcal{H} \nonumber \\
& \phantom{\leq}\, + \tilde{C}_2 \sqrt{1+\epsilon} \sup\limits_{\tau \in [0,t+(D+t_0)]} e^{\delta_\epsilon \tau} \Vert v(\tau) \Vert \nonumber \\
& \leq F_1 \sqrt{V(D+t_0)} 
+ F_2 \Vert x_0 \Vert 
+ F_3 \sup\limits_{\tau \in [0,D+t_0]} e^{\delta_\epsilon \tau} \Vert X(\tau) \Vert_\mathcal{H} \nonumber \\
& \phantom{\leq}\, + F_4 \sup\limits_{\tau \in [0,t+(D+t_0)]} e^{\delta_\epsilon \tau} \Vert v(\tau) \Vert , \label{eq: IDS-ODE stab estimate 1 - intermediate 2}\end{align}
where $F_1 = \tilde{C}_1 E_1 \sqrt{1+\epsilon}$, $F_2 = \tilde{C}_0 + \tilde{C}_1 E_2 \sqrt{1+\epsilon}$, $F_3 = \tilde{C}_1 ( 1 + E_3 ) \sqrt{1+\epsilon}$, and $F_4 = (\tilde{C}_2 + \tilde{C}_1 E_4) \sqrt{1+\epsilon}$. Combining (\ref{eq: IDS-ODE stab estimate 1 - intermediate 1}-\ref{eq: IDS-ODE stab estimate 1 - intermediate 2}) and noting that
\begin{align}
& V(D+t_0) \nonumber \\
& \leq C_9 \Vert X_0 \Vert_\mathcal{H}^2 + C_{10} \underset{\tau \in [0,D+t_0]}{\mathrm{sup}} \Vert d(\tau) \Vert_\mathcal{H}^2 \nonumber \\
& \leq C_9 \Vert X_0 \Vert_\mathcal{H}^2 
+ D_1 C_{10} \underset{\tau \in [0,D+t_0]}{\mathrm{sup}} \Vert x(\tau) \Vert^2 \label{eq: estimate V(D+t0)} \\
& \phantom{\leq} + D_2 C_{10} \underset{\tau \in [0,D+t_0]}{\mathrm{sup}} \Vert X(\tau) \Vert_\mathcal{H}^2
+ D_3 C_{10} \underset{\tau \in [0,D+t_0]}{\mathrm{sup}} \Vert v(\tau) \Vert^2 , \nonumber
\end{align}
we obtain the existence of constants $G_i \geq 0$, independent of the initial condition $(x_0,X_0)$ and the disturbance $v$, such that (\ref{eq: IDS-ODE stab estimate 1}) holds true for all $t \geq D+t_0$. Finally, based on (\ref{eq: estimation of u from V}) and (\ref{eq: conv lemma - applied to V}), we estimate the control input as follows. For all $t \geq 0$,
\begin{align*}
& e^{\delta_\epsilon t} \Vert u(t+(D+t_0)) \Vert \\
& \leq \dfrac{\Vert K \Vert}{\sqrt{C_2(\gamma_1)}} \sqrt{V(t + (D+t_0))} e^{\delta_\epsilon t} \\
& \leq \dfrac{\Vert K \Vert \tilde{C}_0}{\sqrt{C_2(\gamma_1)}} \sqrt{V(D+t_0)} \\
& \phantom{\leq}\, + \Vert K \Vert \sqrt{\dfrac{C_6}{2 \kappa_0 C_2(\gamma_1)}(1+\epsilon)} e^{-\delta_\epsilon t_0} \sup\limits_{\tau \in [t_0,t+(D+t_0)]} e^{\delta_\epsilon \tau} \Vert d(\tau) \Vert_\mathcal{H} .
\end{align*}
Therefore, we infer from (\ref{eq: loop IDS-ODE estimation d}) and (\ref{eq: estimate V(D+t0)}) the existence of constants $H_i$, independent of the initial condition $(x_0,X_0)$ and the disturbance $v$, such that (\ref{eq: IDS-ODE stab estimate 2}) holds true. It concludes the proof of Theorem~\ref{thm: stab IDS-ODE loop}.

\begin{remark}
In the context of Remark~\ref{rem: blow up in finite time}, i.e., when replacing the first point of Assumption~\ref{assum: A4} by the small gain condition (\ref{eq: sufficient small gain condition for no blow up before D+t0}), the reasoning above still applies over the time interval $[0,t_{\max})$ because $t_{\max} > D+t_0$. In this case, estimate (\ref{eq: IDS-ODE stab estimate 1}) holds true for all $t \in [D+t_0,t_{\max})$. As the supremum of the right-hand side of (\ref{eq: IDS-ODE stab estimate 1}) over any time interval $[D+t_0,T]$ of finite length is finite, we deduce that $t_{\max} = + \infty$. Therefore, the conclusion of Theorem~\ref{thm: stab IDS-ODE loop} still holds true.
\end{remark}

\section{Case study}\label{sec: case study}
In this section, $\mathcal{H}$ denotes the $\mathbb{R}$-Hilbert space of square-integrable functions $L^2(0,L)$ endowed with the inner product $\left< f , g \right>_\mathcal{H} = \int_{0}^{L} f g \diff x$. We consider the following coupled system composed of a one-dimensional ODE and a one-dimensional reaction-diffusion equation on $(0,L)$ with delayed Dirichlet boundary controls located at both ends of the domain
\begin{equation*}
\left\{\begin{split}
& \dot{x}(t) =  f_1(x(t),y(t,\cdot),v(t)) \\
& y_t(t,\xi) = a y_{\xi\xi}(t,\xi) + c y(t,\xi) + f_2(x(t),y(t,\cdot),v(t)) \\
& \begin{bmatrix} y(t,0) \\ y(t,1) \end{bmatrix} = u(t-D)   \\
\end{split}\right.
\end{equation*}
where $(t,\xi) \in \mathbb{R}_+ \times (0,L)$, $X(t) = y(t,\cdot) \in \mathcal{H}$, $x(t),v(t) \in \mathbb{R}$, and $u(t) \in \mathbb{R}^2$. The considered coupling functions are given by $f_1(x,X,v) = - a_1 x + \dfrac{b_1}{L} \int_0^L \eta_1 X \diff \xi + c_1 v$ and $f_2(x,X,v) = a_2 x \theta_1 + b_2 \arctan\left( \dfrac{d_2}{L} \int_0^L \eta_2 X \diff \xi \right) \theta_2 + c_2 v \theta_3$ with $a,c,a_i,b_i,c_i,d_i \in \mathbb{R}$, $a,a_1 > 0$, and $\eta_i,\theta_i \in \mathcal{H}$ such that $\Vert \eta_i \Vert_\mathcal{H} = \Vert \theta_i \Vert_\mathcal{H} = 1$.

We define the operator $\mathcal{A}f = a f'' + c f$ over the domain $D(\mathcal{A}) = H^2(0,L)$ and the boundary operator $\mathcal{B}f = (f(0),f(L))$ over the domain $D(\mathcal{B}) = H^1(0,L)$. We introduce the lifting operator $B$ defined for any $(u_1,u_2) \in \mathbb{R}^2$ by $\{B(u_1,u_2)\}(\xi) = u_1 + (u_2 - u_1)\xi/L$ with $\xi \in (0,L)$. We have that the disturbance free operator $\mathcal{A}_0$: 1) generates a $C_0$-semigroup ; 2) is a Riesz-spectral operator with $\lambda_n = c - a n^2 \pi^2 / L^2$ and $\phi_n (\xi) = \psi_n(\xi) = \sqrt{2/L} \sin( n \pi x / L )$, $n \geq 1$. Thus, $(\mathcal{A},\mathcal{B})$ is a boundary control system satisfying Assumptions~\ref{assum: A1} and~\ref{assum: A2}. Furthermore, straightforward computations show that $b_{n,1} = a n \pi \sqrt{2 / L^3}$ and $b_{n,2} = (-1)^{n+1} a n \pi \sqrt{2 / L^3}$. Thus, based on Remark~\ref{rem: A3}, Assumption~\ref{assum: A3} about the Kalman condition is satisfied.

Finally, with the considered coupling functions $f_1$ and $f_2$, Assumption~\ref{assum: A4} holds true with $\tilde{C}_0 = \sqrt{2}$, $\tilde{C}_1 = 2 \vert b_1 \vert / (a_1 L)$, $\tilde{C}_2 = 2 \vert c_1 \vert / a_1$, $D_1 = \vert a_2 \vert$, $D_2 = \vert b_2 d_2 \vert / L$, and $D_3 = \vert c_2 \vert$.

For numerical computations, we take $L = 2 \pi$, $D = 0.1\,\mathrm{s}$, $a = 5$ and $c = 2.5$. Thus, we have one unstable mode with $\lambda_1 = 1.25$ while $\lambda_2 = -2.5$, and $\lambda_3 = - 8.75$. For design purposes, we consider a second order truncated model, i.e., $N_0 = 2$ and $\alpha = 8.75$. Then, the feedback gain matrix $K \in \mathbb{R}^{2 \times 2}$ is  computed based on this truncated model such that the two poles are both placed at $-3$. Following the developments of Section~\ref{sec: Study of the infinite-dimensional closed-loop system}, the degrees of freedom available in the choice of the parameters $\beta \in (0,1)$, $\gamma_1 > C_1/\lambda_m(P)$, and $\gamma_2 > \max \left( \Vert B K \Vert^2/(m_R \lambda_m(P)) , C_5/(1-\beta) \right)$ are used to minimize the value of the constant $C_4 \sqrt{C_6/(2 \kappa_0)}$ involved in the small gain condition (\ref{eq: sufficient small gain condition}). With the \textsc{Matlab} function \texttt{fminsearch}, we obtain with $\beta = 0.4131$, $\gamma_1 = 106.3290$, and $\gamma_2 = 337.1938$ the value $C_4 \sqrt{C_6/(2 \kappa_0)} \approx 8.6260$. Thus, Theorem~\ref{thm: stab IDS-ODE loop} applies when the vector fields $f_1$ and $f_2$ are such that $2 \vert b_1 a_2 \vert / a_1 + \vert b_2 d_2 \vert < L/8.6260 \approx 0.7284$.

\begin{figure}
\centering
\includegraphics[width=3.3in,height=2.25in]{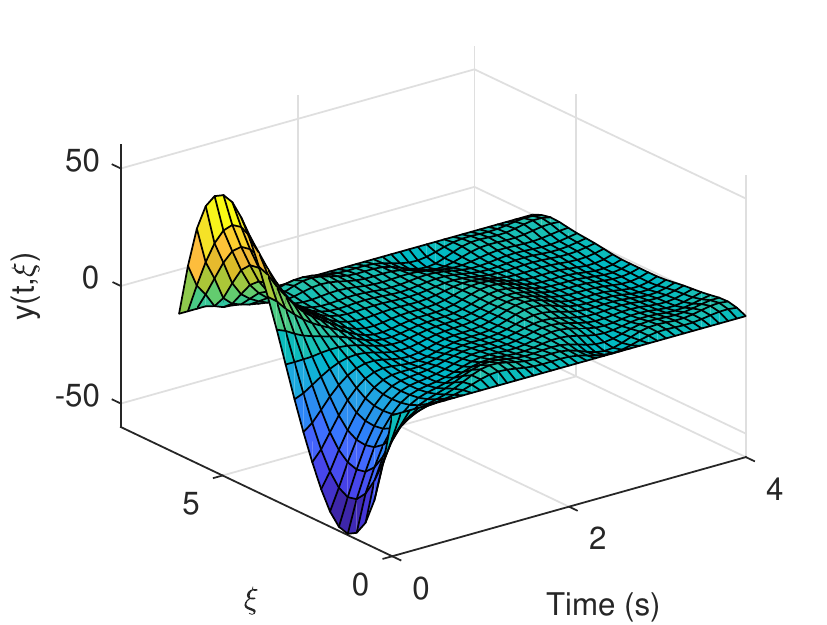}
\caption{Time evolution of the reaction-diffusion part of the closed-loop system}
\label{fig: time evolution X}

\includegraphics[width=3.3in,height=1in]{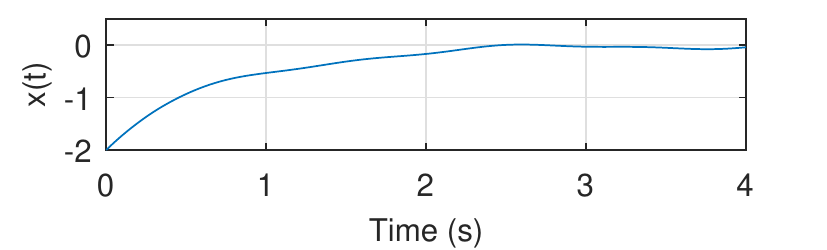}
\caption{Time evolution of the ODE part of the closed-loop system}
\label{fig: time evolution x}

\includegraphics[width=3.3in,height=1in]{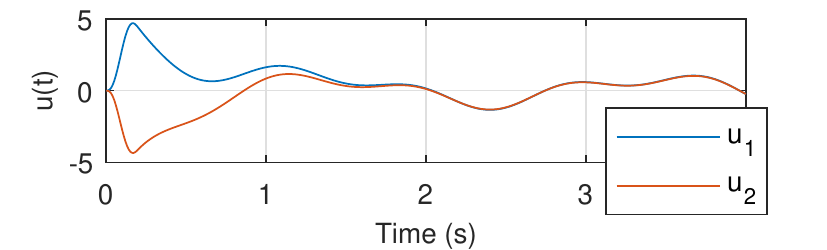}
\caption{Command effort of the closed-loop system}
\label{fig: time evolution u}
\end{figure}

Consequently, we select for numerical simulations $a_1 = 1.5$, $b_1 = 0.5$, $c_1 = 0.2$, $a_2 = 0.7$, $b_2 = 0.55$, $c_2 = 10$, $d_2 = 0.45$, $\eta_1 = \eta_2 = \theta_2 = \sqrt{6 \xi (L - \xi)}/L^{3/2}$, $\theta_1 = \sqrt{2 \xi}/L$, and  $\theta_3 = \sqrt{2 (L-\xi)}/L$. The transition time $t_0$ is set to $t_0 = 0.2\,\mathrm{s}$ while the switching function $\left.\varphi\right\vert_{[0,t_0]}$ is selected as the restriction over $[0,t_0]$ of the unique quintic polynomial function $f$ satisfying $f(0)=f'(0)=f''(0)=f'(t_0)=f''(t_0)=0$ and $f(t_0)=1$. The adopted numerical scheme consists in the discretization of the reaction-diffusion equation using its first 10 modes. The evolution of the closed-loop system is depicted in Figs.~\ref{fig: time evolution X}-\ref{fig: time evolution u} for the initial condition $x_0 = -2$ and $X_0(\xi) = -5\xi(L/2-\xi)(L-\xi)$, and with the external disturbance $v(t) = \sin(2t) \sin(5t)$. The obtained numerical results are compliant with the theoretical predictions.

\section{Conclusion}\label{sec: conclusion}
This paper discussed the feedback stabilization of a class of diagonal Infinite-Dimensional Systems (IDS) with delay boundary control. The proposed approach generalizes a design method formerly reported for a reaction-diffusion equation while proposing a simplification of the boundary control law. The method consists, via a spectral decomposition, in the synthesis of a state-feedback for a finite-dimensional subsystem capturing the unstable dynamics of the plant. Due to the input delay, the design of the control law on the truncated subsystem has been carried out by means of the Artstein transformation. Then, an adequate Lyapunov function has been introduced to assess that the control law designed on the truncated subsystem also ensures the stabilization of the original IDS. Furthermore, it has been shown that this Lyapunov function also allows the assessment of the Input-to-State Stability (ISS) of the closed-loop system with respect to distributed disturbances. Finally, this ISS property has been used to study the stability of the closed-loop IDS when interconnected with an Ordinary Differential Equation (ODE) that also satisfies an ISS property. Specifically, it has been shown that the satisfaction of a certain small gain condition ensures the stability of the IDS-ODE loop for the proposed delayed boundary control law.

% or

%\appendices
%\section{Annex}

\appendix[Regularity and time derivative of an infinite sum]\label{annex1}
Let $\left\{ e_n , \; n \in \mathbb{N}^* \right\}$ by a Hilbert basis of $\mathcal{H}$. Then, as $\left\{ \phi_n , \; n \in \mathbb{N}^* \right\}$ is a Riesz basis with associated biorthogonal set $\left\{ \psi_n , \; n \in \mathbb{N}^* \right\}$, there exists $T \in \mathcal{L}(\mathcal{H})$ such that $T^{-1} \in \mathcal{L}(\mathcal{H})$ and, for all $n \geq 1$, $\phi_n = T e_n$ and $\psi_n = (T^{-1})^* e_n$. Let $A \in \mathcal{C}^1(\mathbb{R}_+;\mathcal{H})$ be given. We obtain that, for all $t \geq 0$,
\begin{align*}
\sum\limits_{k \geq 1} \vert \left< A(t) , \psi_k \right> \vert^2
& = \sum\limits_{k \geq 1} \vert \left< A(t) , (T^{-1})^* e_k \right> \vert^2 \\
& = \sum\limits_{k \geq 1} \vert \left< T^{-1} A(t) , e_k \right> \vert^2 \\
& = \Vert T^{-1} A(t) \Vert_\mathcal{H}^2 \\
& = \left< T^{-1} A(t) , T^{-1} A(t) \right>_\mathcal{H} .
\end{align*}
Thus $\sum\limits_{k \geq 1} \vert \left< A , \psi_k \right> \vert^2 \in \mathcal{C}^1(\mathbb{R}_+;\mathbb{R})$ and we have for all $t \geq 0$,
\begin{align*}
& \dfrac{\mathrm{d}}{\mathrm{d}t} \left[ \dfrac{1}{2} \sum\limits_{k \geq 1} \left\vert \left< A(t) , \psi_k \right>_\mathcal{H} \right\vert^2 \right] \\
& = \operatorname{Re} \left< T^{-1} \dfrac{\mathrm{d}A}{\mathrm{d}t}(t) , T^{-1} A(t) \right>_\mathcal{H} \\
& = \operatorname{Re} \left< T^{-1} \sum\limits_{k \geq 1} \left< \dfrac{\mathrm{d}A}{\mathrm{d}t}(t) , \psi_k \right>_\mathcal{H} \phi_k , T^{-1} \sum\limits_{l \geq 1} \left< A(t) , \psi_l \right>_\mathcal{H} \phi_l \right>_\mathcal{H} \\
& = \sum\limits_{k , l \geq 1} \operatorname{Re} \left\{ \left< \dfrac{\mathrm{d}A}{\mathrm{d}t}(t) , \psi_k \right>_\mathcal{H} \overline{\left< A(t) , \psi_l \right>_\mathcal{H}} \left< T^{-1} \phi_k , T^{-1} \phi_l \right>_\mathcal{H} \right\} \\
& = \sum\limits_{k , l \geq 1} \operatorname{Re} \left\{ \left< \dfrac{\mathrm{d}A}{\mathrm{d}t}(t) , \psi_k \right>_\mathcal{H} \overline{\left< A(t) , \psi_l \right>_\mathcal{H}} \left< e_k , e_l \right>_\mathcal{H} \right\} \\
& = \sum\limits_{k \geq 1} \operatorname{Re} \left\{ \left< \dfrac{\mathrm{d}A}{\mathrm{d}t}(t) , \psi_k \right>_\mathcal{H} \overline{\left< A(t) , \psi_k \right>_\mathcal{H}} \right\} .
\end{align*}
Noting that, for all $k \geq 1$,
\begin{equation*}
\dfrac{\mathrm{d}}{\mathrm{d}t} \left[ \dfrac{1}{2} \left\vert \left< A(t) , \psi_k \right>_\mathcal{H} \right\vert^2 \right]
= \operatorname{Re} \left\{ \left< \dfrac{\mathrm{d}A}{\mathrm{d}t}(t) , \psi_k \right>_\mathcal{H} \overline{\left< A(t) , \psi_k \right>_\mathcal{H}} \right\} ,
\end{equation*}
we deduce that $\sum\limits_{k \geq N_0 + 1} \vert \left< A , \psi_k \right> \vert^2 \in \mathcal{C}^1(\mathbb{R}_+;\mathbb{R})$.

%and we have for all $t \geq 0$, 
%\begin{align*}
%& \dfrac{\mathrm{d}}{\mathrm{d}t} \left[ \dfrac{1}{2} \sum\limits_{k \geq N_0 + 1} \left\vert \left< A(t) , \psi_k \right>_\mathcal{H} \right\vert^2 \right] \\
%& = \sum\limits_{k \geq N_0 + 1} \operatorname{Re} \left\{ \left< \dfrac{\mathrm{d}A}{\mathrm{d}t}(t) , \psi_k \right>_\mathcal{H} \overline{\left< A(t) , \psi_k \right>_\mathcal{H}} \right\} .
%\end{align*}
 
%\section*{Acknowledgment}

\ifCLASSOPTIONcaptionsoff
  \newpage
\fi



\bibliographystyle{IEEEtranS}
\nocite{*}
\bibliography{IEEEabrv,mybibfile}

\end{document}